\documentclass[12pt]{amsart}
\usepackage{frontmatter}

\usepackage[numbers,comma,square,sort&compress]{natbib} 
\numberwithin{equation}{section}

\title[Micropterons, Nanopterons and Solitary Waves for Diatomic FPUT]{Micropterons, Nanopterons and Solitary Wave Solutions to the Diatomic Fermi-Pasta-Ulam-Tsingou Problem}

\author{Timothy E. Faver}
\address{Mathematical Institute, Universiteit Leiden, P.O. Box 9512, 2300 RA Leiden, The Netherlands, {\tt{t.e.faver@math.leidenuniv.nl}}}

\author{Hermen Jan Hupkes}
\address{Mathematical Institute, Universiteit Leiden, P.O. Box 9512, 2300 RA Leiden, The Netherlands, {\tt{hhupkes@math.leidenuniv.nl}}}



\date{\today}

\begin{document}

\begin{abstract}
We use a specialized boundary-value problem solver for mixed-type functional differential equations to numerically examine the landscape of traveling wave solutions to the diatomic Fermi-Pasta-Ulam-Tsingou (FPUT) problem. 
By using a continuation approach, we are able to uncover the relationship between the branches of micropterons and nanopterons that have been rigorously constructed recently
in various limiting regimes. 
We show that the associated surfaces are connected together in a nontrivial fashion and illustrate the key role that solitary waves play in the branch points. 
Finally, we numerically show that the diatomic solitary waves are stable under the full dynamics of the FPUT system.
\end{abstract}

\maketitle

\section{Introduction}
In this paper we study the Fermi-Pasta-Ulam-Tsingou (FPUT) equation
\cite{fput-original,dauxois}
\begin{equation}\label{newton}
m_j \ddot{x}_j
= F(x_{j+1} - x_j) -F(x_{j} - x_{j-1}),
\qquad \qquad j \in \mathbb{Z}
\end{equation}
in the diatomic regime
\begin{equation}\label{masses}
m_j = 
\begin{cases}
1, &j \text{ is odd}, \\
m, &j \text{ is even},
\end{cases}
\end{equation}
using the quadratic spring force $F(r) = r + r^2$. 
In particular, we numerically investigate several branches of diatomic traveling waves that have recently been shown to exist for this system. 
Using a continuation approach, we track these branches outside of the parameter regimes where they have been rigorously constructed.  
This allows us to shed light on the intricate structure of the broader landscape of such solutions.

\subsection{Propagation in discrete media}
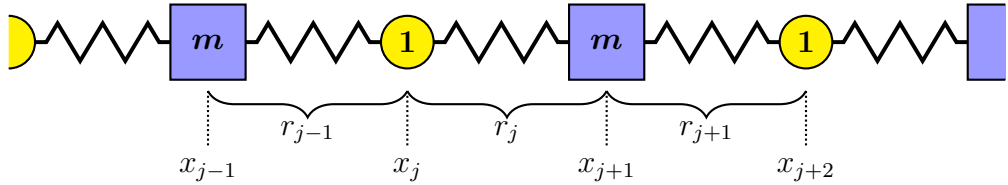
\begin{figure}
\[
\begin{tikzpicture}

%%------------------------------------------------------------------------------------------------------------------------------------------------------------------------------------------------------------%%
\def\arith#1#2#3{
% Compute the product of the sum #3*(#1+#2).
% Tikz arithmetic does not like distributing!
#3*#1+#3*#2
};

%%------------------------------------------------------------------------------------------------------------------------------------------------------------------------------------------------------------%%
\def\squaremass#1#2#3{
% Draws a square with the following coordinates:
% #1 = upper left X-coordinate of square
% #2 = half of side length of square
% #3 = label inside square
\draw[fill=blue, opacity=.4, draw opacity=1,thick]
(#1,#2) rectangle (#1+2*#2,-#2);
\node at (#1+#2,0){$\boldsymbol{#3}$};
};

%%------------------------------------------------------------------------------------------------------------------------------------------------------------------------------------------------------------%%
\def\circlemass#1#2#3{
% Draws a circle with the following parameters:
% #1 = X-coordinate of center (Y-coordinate = 0)
% #2 = radius
% #3 = label inside circle
\draw[fill=yellow,thick] (#1,0)node{$\boldsymbol{#3}$} circle(#2);
};

%%------------------------------------------------------------------------------------------------------------------------------------------------------------------------------------------------------------%%
\def\spring#1#2#3#4{
% Draws a pointy spring with three ``upper points'' and two ``lower points'' and ``connectors'' on the left and right.
% #1 = X-coordinate of left endpoint of spring
% #2 = connector length
% #3 = ``wavelength'' of spring = distance between a ``pointy peak'' and the next X-axis intersection
% #4 = height of spring
\draw[line width = 1.5pt] (#1,0)
--(#1 + #2,0)
--(#1 + #2 + #3, #4)
--(#1 + #2 + 3*#3,-#4)
--(#1 + #2 + 5*#3,#4)
--(#1 + #2 + 7*#3,-#4)
--(#1 + #2 + 9*#3,#4)
--(#1 + #2 + 10*#3,0)
--(#1 + 2*#2 + 10*#3,0);
};

%%------------------------------------------------------------------------------------------------------------------------------------------------------------------------------------------------------------%%
\def\d{.15}; % horizontal spring connector length
\def\mh{1}; % side length of square mass
\def\sh{.25}; % height of spring 
\def\swl{.15}; % wave length of spring
\def\r{.35}; %radius of circular mass
\def \sl#1{\arith{2*\d}{10*\swl}{#1}}; % spring length, multiplied by the factor #1

%%------------------------------------------------------------------------------------------------------------------------------------------------------------------------------------------------------------%%
\def\labeldown#1#2{
% Draws a densely dotted vertical line with X-coordinate at #1 and Y-coordinate fixed.
% Puts the label #2 at the bottom of the line
\draw[densely dotted,thick] (#1,-\mh/2-\d)--(#1,-\mh/2-\d-.75)node[below]{#2}
};

%%------------------------------------------------------------------------------------------------------------------------------------------------------------------------------------------------------------%%
\def\brace#1#2{
% Draws a curly brace with left endpoint X-coordinate equal to #1 and all other coordinates fixed.
% Puts the label #2 at the middle.
\draw[decoration={brace, amplitude = 10pt,mirror},decorate,thick] 
(#1,-\mh/2-\d)--node[midway,below=7pt]{#2}(#1+\sl{1}+\mh/2+\r,-\mh/2-\d)
}

% LEFT HALF CIRCULAR MASS

\fill[yellow] (-\r,\r) arc(90:-90:\r)--cycle;
\draw[thick] (-\r,\r) arc(90:-90:\r);

% FIRST SPRING
\spring{0}{\d}{\swl}{\sh};

% FIRST SQUARE MASS

\squaremass{\sl{1}}{\mh/2}{m};

% SECOND SPRING

\spring{\sl{1}+\mh}{\d}{\swl}{\sh};

% FIRST CIRCULAR MASS

\circlemass{\sl{2}+\mh+\r}{\r}{1};

% THIRD SPRING

\spring{\sl{2}+\mh+2*\r}{\d}{\swl}{\sh};

% SECOND SQUARE MASS

\squaremass{\sl{3}+\mh+2*\r}{\mh/2}{m};

% FOURTH SPRING

\spring{\sl{3}+2*\mh+2*\r}{\d}{\swl}{\sh};

% SECOND CIRCULAR MASS

\circlemass{\sl{4}+2*\mh+3*\r}{\r}{1};

% FIFTH SPRING

\spring{\sl{4} + 2*\mh+4*\r}{\d}{\swl}{\sh};

% RIGHT HALF MASS

\fill[blue,opacity=.4]
(\sl{5}+2*\mh+4*\r,\mh/2) rectangle (\sl{5}+2*\mh+4*\r+\mh/2,-\mh/2);

\draw[thick] (\sl{5}+2*\mh++4*\r+\mh/2,\mh/2)
--(\sl{5}+2*\mh+4*\r,\mh/2)
--(\sl{5}+2*\mh+4*\r,-\mh/2)
--(\sl{5}+2*\mh+4*\r+\mh/2,-\mh/2);

% LABELS

\labeldown{\sl{1}+\mh/2}{$x_{j-1}$};
\labeldown{\sl{2}+\mh+\r}{$x_j$};
\labeldown{\sl{3}+3*\mh/2+2*\r}{$x_{j+1}$};
\labeldown{\sl{4}+2*\mh + 3*\r}{$x_{j+2}$};

\brace{\sl{1}+\mh/2}{$r_{j-1}$};
\brace{\sl{2}+\mh+\r}{$r_j$};
\brace{\sl{3}+3*\mh/2+2*\r}{$r_{j+1}$};

\end{tikzpicture}
\]
\caption{The diatomic FPUT lattice, featuring alternating masses $m$ and 1 connected by identical springs.}
\label{fig-lattice_fig}
\end{figure}

By now, the FPUT system \eqref{newton} is well-established as a classic prototype of the propagation of disturbances through spatially discrete systems, such as granular media, artificial metamaterials, DNA strands, and electrical transmission lines \cite{brill, kev}. 
In essence, it models an infinite, one-dimensional chain of particles that can only move horizontally and are connected to their nearest left and right neighbors by nonlinear springs. 
In the relative displacement coordinates
\begin{equation}
r_j 
=x_{j+1}  - x_j,
\end{equation}
these springs transmit a force $F(r_j)$ between the particles at sites $j$ and $j+1$. Applying Newton's law, this leads naturally to the evolution \eqref{newton} for the position variables $x_j$; see Fig.\@ \ref{fig-lattice_fig}. 
We do note that various other studies incorporate higher-order terms into the force $F(r)$ \cite{friesecke-pego1, herrmann-matthies-asymptotic, pankov}, but these can typically be handled using refined perturbative techniques \cite{faver-spring-dimer} provided that the displacements from equilibria remain sufficiently small. 

Traveling waves have played a fundamental role in the analysis of \eqref{newton} and other spatially discrete systems \cite{kev,MPB,CGW08,HJHFZHNGM,VL28}. 
In the current diatomic setting, such solutions take the form
\begin{equation}\label{tw ansatz}
r_j(t) 
= \begin{cases}
\rbar_o(j + ct), &j \text{ is odd}, \\
\rbar_e(j + ct), &j \text{ is even},
\end{cases}
\end{equation}
which means that the speed $c$ and the pair of waveprofiles $(\rbar_o, \rbar_e)$ must satisfy the two-component mixed-type functional differential equation\footnote{Such equations are also referred to as advance-delay differential equations.} 
(MFDE)
\begin{equation}
\label{tw eqns}
\begin{cases}
c^2\rbar_o''(\xi)
= \frac{1}{m}F\big(\rbar_e(\xi +1)\big) - \left(1+\frac{1}{m}\right) F\big(\rbar_o(\xi)\big) + F\big(\rbar_e(\xi -1)\big)
\\
c^2\rbar_e''(\xi) 
= F\big(\rbar_o(\xi +1)\big) - \left(1+\frac{1}{m}\right)F\big(\rbar_e(\xi) \big) + \frac{1}{m}F\big(\rbar_o(\xi -1)\big).
\end{cases}
\end{equation}
The argument shifts in this system prevent the use of powerful ODE-techniques such as phase-plane analysis. 
This causes many technical complications in the analysis of spatially discrete systems; see, e.g., \cite{hupkes2019traveling} for an overview of the machinery that has been developed for MFDEs.

\subsection{Formal limits}
The complexity of  \eqref{tw eqns} can be reduced considerably by taking various (formal) limits. 
We briefly sketch three procedures of this type that play an important role in this paper.

\subsubsection{Equal mass limit}\label{equal mass limit intro}
Taking $m=1$, the diatomic lattice becomes monatomic, and one can set $\rbar_o = \rbar_e = \phi$ to arrive at the scalar MFDE
\begin{equation}\label{mono:fput:tw:eqn}
c^2 \phi''(\xi) 
= F\big(\phi(\xi+1)\big) - 2 F\big(\phi(\xi)\big) + F\big(\phi(\xi - 1) \big) .
\end{equation}
This is a classical problem that was analyzed in detail by Friesecke in combination with Wattis \cite{friesecke-wattis} and Pego \cite{friesecke-pego1, friesecke-pego2, friesecke-pego3, friesecke-pego4}, who established that there exists a smooth branch of nontrivial solitary waves $(c,\phi_c)$ that are even, exponentially localized, and stable. 
In particular,
for small $\epsilon > 0$ one can write
\begin{equation}\label{fr-p speed of sound}
\cep = 1 + \frac{\ep^2}{24},
\end{equation}
and show that the associated profiles $\phi_{c_{\epsilon}}$ satisfy the limiting behavior \cite{friesecke-pego1}
\begin{equation}\label{fp lw limit}
\norm{\frac{8}{\ep^2} \phi_{\cep}(2\ep^{-1} \cdot ) - \sech^2(\cdot)}_{H^1} 
= \O(\ep^2).
\end{equation}
More recently, Herrmann and Matthies \cite{herrmann-matthies-asymptotic, herrmann-matthies-uniqueness, herrmann-matthies-stability} considered the ``high-energy'' limit $c \gg 1$, using a different Lennard-Jones-type potential for the springs  that is analytic at $r=0$ but singular at $r=1$.

\subsubsection{Small mass limit}\label{small mass section intro}
Multiplying \eqref{tw eqns} by $m$ and formally setting $m = 0$, the injectivity of $F$ on $[0,\infty)$ yields the identification
\begin{equation}
\label{eq:id:small:mass}
\rbar_o(\xi) 
= \rbar_e(\xi+1).
\end{equation}
Physically, this means that the mass-less particles are fixed halfway between the heavier ones, corresponding with the intuition developed in  \cite{hoffman-wright,pelinovsky-schneider}.
Upon setting 
\begin{equation}\label{small mass big R}
\varphi(\xi) 
=\frac{1}{2} \rbar_o(2\xi) + \frac{1}{2}\rbar_e(2\xi+1)
\end{equation}
and adding the first line of \eqref{tw eqns} to a shifted version of the second line, the identification \eqref{eq:id:small:mass} readily reveals the MFDE
\begin{equation}\label{small mass MFDE}
c^2 \varphi''(\xi) 
= 2\big[F\big(\varphi(\xi + 1)\big)- 2 F\big( \varphi(\xi)\big) + F\big( \varphi(\xi - 1)\big)\big].
\end{equation}
Comparing this with \eqref{mono:fput:tw:eqn} immediately shows that one may write $\varphi_{c} = \phi_{\sqrt{2} c}$, which means that solitary waves can be expected for $c \gtrsim \sqrt{2}$.
This extra scaling corresponds to the notion that the effective limiting monatomic lattice has double the spring length of the original diatomic lattice.

\subsubsection{Long wave limit}\label{lw limit sol section}
Here we fix $m \ne 1$ and make the classical long wave scaling \cite{schneider-wayne-water-wave}
\begin{equation}
\rbar_o(\xi) = \ep^2\theta_o(\ep\xi;\ep)
\quadword{and}
\rbar_e(\xi) = \ep^2\theta_e(\ep\xi;\ep).
\end{equation}
Upon making the further perturbation ansatz
\begin{equation}
\theta_o(X;\ep)
= \sum_{k=0}^3 \ep^k u_k(X)
\quadword{and}
\theta_e(X;\ep)
= \sum_{k=0}^3 \ep^k v_k(X),
\end{equation}
one can subsequently solve the traveling wave equations \eqref{tw eqns} formally to $\O(\ep^4)$ by taking
\begin{equation}
u_1 = v_1 = \Phi_m, 
\qquad 
\qquad
c = c_m^{(s)}  + \O(\ep^2)
\end{equation}
and defining the other $u_k$ and $v_k$ in more complicated terms of the profile $\Phi_m$.
Here the speed of sound is
\begin{equation}\label{lw speed of sound}
c^{(s)}_m
:= \sqrt{\frac{2}{1+m}}
\end{equation}
and the solitary wave profile $\Phi_m$ satisfies the KdV traveling wave equation
\begin{equation}
a_m\Phi'' - \Phi + b_m\Phi^2 = 0,
\end{equation}
where $a_m$ and $b_m$ are (complicated) $m$-dependent coefficients, see \cite[Eq.\@ (3.2)]{faver-wright}.
That is,
\begin{equation}
\label{eq:int:def:Phi:m}
\Phi_m(X)
:= \frac{3}{2b_m}\sech^2\left(\frac{X}{2\sqrt{a_m}}\right).
\end{equation}
We note that $c^{(s)}_m$ reduces to the critical values found above for $m=1$ and $m=0$.

In a certain sense, this procedure can be seen as a specialization of the techniques developed in \cite{schneider-wayne} and \cite{gmwz}. 
Here the authors consider monatomic respectively polyatomic\footnote{In this case both the masses and the spring forces in \eqref{newton} are allowed to vary periodically.} FPUT systems and derive a set of KdV PDEs to approximate the evolution of suitably scaled initial conditions over algebraically long time-scales.

\subsection{Rigorous results}
\begin{figure}[t]
\[
\begin{tikzpicture}[thick]

\draw[white] (-1,7) rectangle (10,-2);

% AXES

\draw[<->] (0,5.5)node[right]{$m$}--(0,0)--(6.5,0)node[right]{$c$};

\draw (0,4.5)--(.25,4.5)node[right]{1};
\draw (1,0)--(1,-.25)node[below]{1};

\draw (1.85,0)--(1.85,-.25)node[below]{$\sqrt{2}$};

% LONG WAVE

\begin{scope}

\clip (0,4.5)--(0,0)--(1.85,0) to[bend left=17] (4.5+.85,1)--(4.5+.85,4.5)--cycle;

\fill[yellow,opacity=.4] (1,4.5) to[bend right = 10] (1.85,0)--(4,0) to[bend left = 20] (1,4.5);

\draw[blue, line width = 3pt] (1,4.5) to[bend right=10] (1.85,0);
\draw (1,4.5) to[bend right=20] (4,0);

\end{scope}

\begin{scope}

\clip (1,4.5) to[bend right = 10] (1.85,0)--(4,0) to[bend left = 20] (1,4.5);
\draw[densely dotted] (1,.85)--(4,.85);
\end{scope}

\node[circle,draw=blue,outer sep=0pt,inner sep = 0pt,fill=white] at (1.5,.85){$\Phi_{\!m}$};

\draw[<-] (1.6,2)--(2.6,2)node[right]{long wave limit};

\draw[<-] (1.2,2.5)--(2.25,2.5)node[right]{$\sech^2$-type KdV solitary waves ($c=c_m^{(s)}$)};

% SMALL MASS

\fill[red,opacity=.4] (1.85,0) to[bend left=17] (4.5+.85,1)--(4.5+.85,0)--cycle;

\draw (1.85,0) to[bend left=17] (4.5+.85,1);

\draw[blue,line width = 3pt] (1.85-.015,0)--(4.5+.85,0);
\draw[<-] (4+.85,.75)--(4.7171+.85+.2,1.75-.2)node[right]{small mass};
\node[right] at (4.7171+.85+.2,1.4-.2){limit};
\draw[<-] (3.75+.85,-.1)--(3.75+.85,-1.01)node[below]{monatomic solitary waves};

\begin{scope}
\clip (1+.85,0) to[bend left=17] (4.5+.85,1)--(4.5+.85,0)--cycle;

\draw[densely dotted] (3+.85,1)--(3+.85,0);
\end{scope}

\node[circle,draw=blue,outer sep=0pt,inner sep = 0pt,fill=white] at (3+.85,0){$\varphi_c$};

\draw[dashed] (4.5+.85,0)--(4.5+.85,1);

%% EQUAL MASS

\draw[thick] (1,4.55) to[bend left=8] (4.5,5);
\draw[thick] (1,4.45) to[bend right=8] (4.5,4);

\fill[leidenblue,opacity=.4] (1,4.55) to[bend left=8] (4.5,5)--(4.5,4) to[bend left=8] (1,4.45);
\fill[leidenblue,opacity=.2] (4.5,5) rectangle (6.5,4);

\draw[blue,line width = 3pt] (.985,4.5)--(4.5,4.5);
\draw[blue,dashed,->, line width = 1.5pt] (4.5,4.5)--(6.5,4.5);

\draw[dashed] (4.5,5) rectangle (6.5,4);
\draw[->] (2.75,5.75)node[above]{equal mass limit {($c \gtrsim 1$)}}--(2.75,4.75);
\draw[<-] (4.75,4.6)--(4.75+.7071,5.6)node[right]{monatomic solitary};
\node[right] at (4.75+.7071,5.25){waves {(hypothesis)}};
\draw[<-] (3.5,4.4) to (3.5,3.4);
\node[below] at (5,3.4){monatomic solitary waves};
\draw[<-] (5.8,4.2)--(6.8,4.2)node[right]{equal mass};
\node[right] at (6.8,3.82){limit {($c \gg 1$)}};

\begin{scope}
\clip (1,4.55) to[bend left=8] (4.5,5)--(4.5,4) to[bend left=8] (1,4.45);
\draw[densely dotted] (3,4)--(3,5);
\end{scope}

\node[circle,draw=blue,outer sep=.25pt,inner sep = .25pt,fill=white] at (3,4.5){$\phi_c$};

\end{tikzpicture}
\]
\caption{Bands of rigorously constructed traveling wave solutions for the long wave (yellow), small mass (red), and equal mass (blue) diatomic problems.
In general, it is a technical artifact of the various existence proofs in \cite{faver-wright, hoffman-wright, faver-hupkes} that the bands collapse as the fixed parameter approaches its own critical value.}
\label{fig-all_nonlocals}
\end{figure}
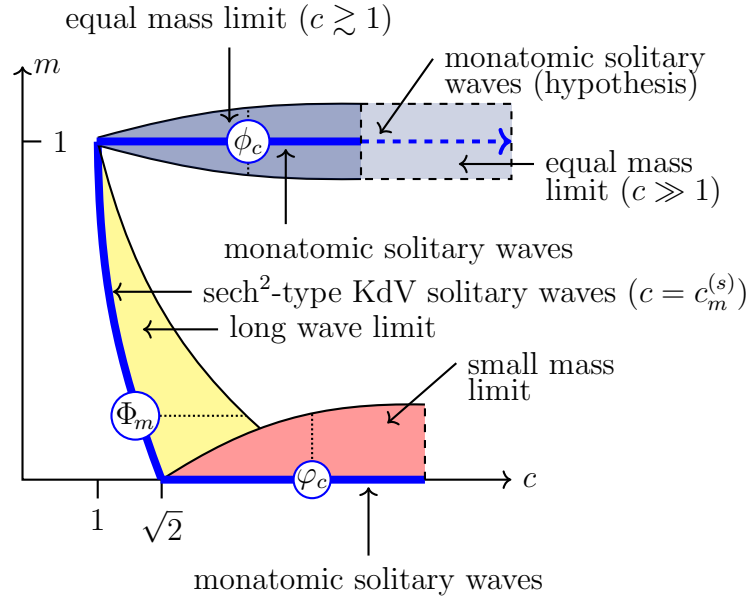

The main focus of the recent papers \cite{faver-wright, hoffman-wright, faver-hupkes} has been to rigorously establish the presence of solutions to the two-component traveling wave problem \eqref{tw eqns} in the neighborhood of the (formal) limiting solutions discussed above. 
The parameter regimes that have now been treated are depicted in Fig.\@ \ref{fig-all_nonlocals}, which we reproduce from 
\cite[Fig.\@ 2]{faver-hupkes}. 
Each of these regimes has its distinctive features and requires specialized tools and techniques, which we briefly discuss below and in 
{\S}\ref{sec:bck}.

\subsubsection{Ripples}
The common theme in the approaches of \cite{faver-wright, hoffman-wright, faver-hupkes} is that one has to give up on the exponential localization of the wave profiles. 
Stated more informally, the traveling wave equations \eqref{tw eqns} are not generically expected to admit solitary wave solutions.
One can interpret this as a manifestation of the purely imaginary spatial eigenvalues associated to the linearization of  \eqref{tw eqns} around the zero background state. 
In contrast to the monatomic setting, this contribution to the essential spectrum cannot be fully removed by applying exponential weight functions. 
Indeed, the resulting linearization typically has codimension one; see {\S}\ref{sec:sub:lin:phi:c}.

To fill the resulting gap, one needs to follow Beale's key insight \cite{beale} and incorporate the background sinusoidal periodic solutions associated to the eigenvalues mentioned above. 
At the nonlinear level, this results in an asymptotic ``ripple'' in the traveling wave profiles at spatial infinity, which destroys their exponential localization.

Quantifying the size of this ripple in terms of the relevant small parameter is an interesting aspect, both from a theoretical and a practical perspective. 
Indeed, the simulations by Giardetti, Shapiro, Windle, and Wright \cite{co-ops} suggest that the size of the ripple is directly related to the energy loss that the ``core'' of the wave experiences as it propagates through the lattice.
We discuss this issue in detail in {\S}\ref{sec:full:fpu}. 
For now, we simply remark that the (meta-)stability of the diatomic waves constructed in \cite{faver-wright, hoffman-wright, faver-hupkes} is a delicate open question.

\subsubsection{Nanopterons}
The results\footnote{It is also possible to find long wave nanopterons in the spring dimer lattice \cite{faver-spring-dimer}.  This is an FPUT lattice in which the masses are identical but the spring forces alternate.} of Faver and Wright in \cite{faver-wright} show that for fixed $m \in (0,1)$  there exist long wave solutions to \eqref{tw eqns} for $c \in (c_m^{(s)}, c_m^{(s)}+\delta_m)$ where $c_m^{(s)}$ is the speed of sound from \eqref{lw speed of sound}. 
Moreover, $\delta_m \to 0$ as $m \to 1^-$.
It is unclear what the behavior of $\delta_m$ is for $m \to 0$, since the relevant parameter in \cite{faver-wright} is in fact $1/m$ which diverges.
The amplitude of the ripple turns out to be small beyond all orders in $c - c_m^{(s)}$, owing to the fact the underlying perturbation from the limiting (scaled) profile $\Phi_m$ is singular.
Following Boyd's \cite{boyd} terminology, we refer to the resulting waves as nanopterons; see Fig.\@ \ref{fig-nano} for a contrast with the solitary wave.

Similar results due to Hoffman and Wright hold for the small mass regime \cite{hoffman-wright}, but here one fixes $c \gtrsim \sqrt{2}$ and the small parameter is $m > 0$.
The difference is that now an underlying solvability condition forces a countable set of mass ratios to be excluded from the analysis. 
It is unclear whether or not (solitary) waves exist at these ratios, which aggregate at zero.

\subsubsection{Micropterons}
We considered the equal mass regime in \cite{faver-hupkes} and established a set of technical conditions under which an (arbitrary) solution $(\phi,c)$ to the monatomic problem \eqref{mono:fput:tw:eqn} can be extended into the setting $m \approx 1$. 
A major difference with the previous settings is that the underlying perturbation problem is regular, which allowed us to provide an explicit integral expression to characterize the $\O(m -1)$ behavior of the ripple amplitude; see {\S}\ref{sec: micropt ampl}.

\begin{figure}
\begin{subfigure}[t]{.3\textwidth}
\centering
\begin{tikzpicture}[scale=.5]

\draw[thick,<->] (-4.5,0)--(4.5,0)node[right]{$x$};
\draw[thick,<->] (0,-.5)--(0,2.5);

\def\a{20};
\def\b{1000};
\def\c{20};
\def\d{1.75};
\def\e{75};
\def\f{1.25};

\draw[ultra thick,blue] plot[domain=-4.25:4.25,smooth,samples=500] (\x,{1.8/(1+3*(\x)^2)});

\end{tikzpicture} 
\caption{A solitary wave}
\label{fig-sol}
\end{subfigure}
\hspace{\fill}
\begin{subfigure}[t]{.3\textwidth}
\centering
\begin{tikzpicture}[scale=.5]
% The core --- drawing

\draw[thick,<->] (-4.5,0)--(4.5,0)node[right]{$x$};
\draw[thick,<->] (0,-.5)--(0,2.5);%node[anchor = west]{$\Sigma_j(X) + \Phi_j(X)$};
%\draw[blue,ultra thick,<->] plot[domain=-4:4, samples=100] (\x,{2.75/(exp(\x)+exp(-\x))});
\draw[ultra thick,blue] plot[domain=-4.25:4.25,smooth,samples=500] (\x,{1.8/(1+3*(\x)^2)});

% The ripple inset --- drawing

\draw[densely dotted,line width = .75pt] (3,{2/(exp(6.25)+exp(-6.25))}) circle(.5);
\draw[densely dotted, line width = .75pt] (3-.5,{2/(exp(5.75)+exp(-5.75))}) -- (5.0253-3.25,2.25);
\draw[densely dotted, line width = .75pt] (3+.5,{2/(exp(6.75)+exp(-6.75))}) -- (7.4747-3.25,2.25);
\draw[densely dotted, line width = .75pt] (3,2.5) circle(1.25);

\begin{scope}
\clip (3,2.5) circle(1.25);
\node at (3-.25,2.5){\begin{tikzpicture}[scale=.5]
\draw[ultra thick,blue] plot[domain=-3:3,samples=100,smooth] (\x,{-.3*cos(10*\x r)});
\end{tikzpicture}};
\end{scope}

\end{tikzpicture}
\caption{A nanopteron}
\label{fig-nano}
\end{subfigure}
\hspace{\fill}
\begin{subfigure}[t]{.3\textwidth}
\centering
\begin{tikzpicture}[scale=.5]

\draw[thick,<->] (-4.5,0)--(4.5,0)node[right]{$x$};
\draw[thick,<->] (0,-.5)--(0,2.5);

\def\a{35};
\def\b{1000};
\def\c{20};
\def\d{1.75};
\def\e{75};
\def\f{1.25};

\draw[ultra thick,blue] plot[domain=-4.25:4.25,smooth,samples=500] (\x,
{
\a*(1-exp(-abs(\x)/\b))*cos(deg(\c*\x))
+ 1.8/(1+3*(\x)^2)
%+ \d*exp(-abs(\x)/\e)/(1+(\f*\x)^2)
});

\end{tikzpicture}
\caption{A micropteron}
\label{fig-micro}
\end{subfigure}

\caption{Schematic representation of the three types of traveling waves featured in this paper.}
\label{fig-all_profiles}
\end{figure}
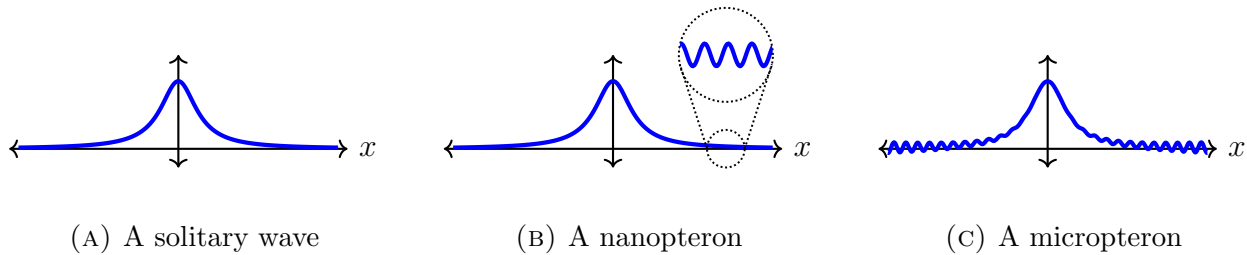

We were able to verify these technical conditions for the waves \eqref{fp lw limit} by performing a careful expansion in $\ep$ and examining the leading order terms. 
However, we strongly suspect that the $\O(m -1)$ coefficient for the ripple amplitude is small beyond all orders in $\ep$, which suggests that our expansions cannot rule out that this coefficient vanishes.
Nevertheless, we followed the terminology of Boyd \cite{boyd} and speculatively referred to our constructed solutions as micropterons; see Fig.\@ \ref{fig-micro} for a contrast with the nanopteron.

\subsubsection{Solitary waves}
Although the ripple amplitude is a crucial variable to close the fixed-point arguments in \cite{faver-wright,hoffman-wright,faver-hupkes}, it could still potentially vanish at certain $(c,m)$ pairs.
The resulting solitary waves can be seen as a lossless mechanism to transport finite-energy states over arbitrary distances. 
As a consequence, solutions of this type play an important role in many applications and have been extensively studied in many different settings \cite{sandstede2002stability,kev, hochstrasser1989energy,jones1991construction}.
It has been conjectured that solitary waves do exist in the diatomic lattice for a countable, discrete set of mass ratios that accumulate at 0; see \cite{faver-wright, VSWP, lustri-porter}.
We discuss these conjectures in greater detail in {\S}\ref{related models section}.

\subsection{Main results}
Besides the justification of our ``micropteron'' terminology from \cite{faver-hupkes}, the main goal of this paper is to numerically examine the full region between the three limiting curves in Fig.\@ \ref{fig-all_nonlocals}. 
In particular, we shed light on the relation between the three types of diatomic waves discussed above by extending them beyond the parameter regimes that were rigorously analyzed in \cite{faver-wright,hoffman-wright,faver-hupkes}. 
Hopefully this will set the stage for further analytical work in this intriguing but challenging area; see, e.g., \cite[{\S}7]{faver-wright} for a short discussion of the technical obstructions.

Our main technical contribution is that we map out the two-parameter surfaces of micropterons and nanopterons that emerge from the (formal) limiting monatomic profiles
at $m=1$, respectively $m=0$ (the horizontal solid curves in Fig.\@ \ref{fig-all_nonlocals}). 
This is achieved by numerically solving the MFDE \eqref{tw eqns} for a large number of $(c,m)$ pairs, using a continuation approach to provide suitable initial conditions. 
In addition, we numerically evaluate the leading-order coefficient for the micropteron ripple-amplitude and show that it does not vanish.

The main conclusion is that the micropteron and nanopteron surfaces are in fact connected, albeit in a highly nontrivial fashion featuring holes, folds, and twists. 
For example, if one fixes $c$ and looks at the associated one-dimensional cross-section, the nanopteron and micropteron curves are typically disconnected. 
Exceptions occur at isolated values of $c$ and appear to be closely related to the occurrence of double roots as one tracks the ripple amplitude over the curve. 

Indeed, as a byproduct of our analysis we uncover several branches of solitary waves. 
These occur in very narrow bands of $m$-values that compress as $m \to 0^+$, which hence provides a numerical confirmation of the conjectures from \cite{VSWP, lustri-porter, lustri} that we discuss below.
In order to explore the stability of these solitary waves, we use (an approximation of) one of them as an initial condition for the dynamics of the full FPUT system \eqref{newton}. 
By comparing the resulting behavior with the computations from \cite{co-ops}, we are able to provide strong numerical evidence to suggest that these waves provide stable and robust mechanisms for energy transport.

\subsection{Related models}\label{related models section}
In order to place our results in perspective, we briefly discuss several prior numerical and theoretical studies that are closely related to our setting. 
These focus on several lattice models that are qualitatively similar to the FPUT system \eqref{newton}.

\subsubsection{The diatomic Toda lattice}
In this setting the spring force in \eqref{newton} is chosen to be the Toda force $F(r) = 1-e^{-r}$ \cite{toda}.
The extra higher-order terms added to our quadratic expression ensure that the problem is fully integrable.
Vainchtein, Starosvetsky, Wright, and Perline \cite{VSWP} study this lattice in the small mass limit.
They use a multiscale asymptotic analysis to deduce that small-amplitude traveling waves in the lattice can have oscillations at infinity unless a certain function $\iota(m)$ of the mass ratio $m$ vanishes.
They then compute the roots of $\iota(m)$ numerically and conjecture that $\iota(m)$ vanishes for a countable number of mass ratios accumulating at 0.
This resembles the FPUT small mass limit in \cite{hoffman-wright}, which could not construct nanopterons at a similar set of mass ratios.
Lustri and Porter \cite{lustri-porter} also work with the diatomic Toda lattice in the small mass limit and use exponential asymptotics  to capture the leading order asymptotics of exponentially small terms in solutions.
They too calculate a countable number of mass ratios for which these terms vanish and only solitary waves should exist.
Lustri \cite{lustri} uses the same techniques for the diatomic FPUT lattice with our quadratic spring force and makes the same conjecture.

\subsubsection{Mass-in-mass lattices}
The mass-in-mass (MiM) lattice is a monatomic lattice of ``beads'' that are ``hollow'' and contain an additional resonator particle \cite{cpkd, kev}.
The spring force connecting the beads is typically Hertzian, which is not smooth, unlike the FPUT and Toda forces.
Various numerical studies predict the formation of nanopterons in Hertzian lattices \cite{kev-vain-et-al, ksx, woodpile, vorotnikov-et-al, xu}.
Conversely, for a countable number of bead-resonator mass ratios accumulating at 0, Kevrekedis, Stefanov, and Xu \cite{ksx} prove the existence of solitary wave solutions for the Hertzian MiM lattice.
Subsequently, Faver, Goodman, and Wright \cite{faver-goodman-wright} have found solitary waves at those same mass ratios when the MiM lattice has the quadratic FPUT spring force.
For mass ratios small but away from this countable set, Faver \cite{faver-mim-nanopteron} has shown the existence of nanopterons, similar to the small mass FPUT limit of \cite{hoffman-wright}.

\subsubsection{$1\!:\!N$ dimers}
The ``$1\!:\!N$ dimer'' is a polyatomic lattice in which one ``heavy'' mass alternates with $N$ ``light'' masses.
Jayaprakash, Vakakis, and Starosvetsky numerically observe a decreasing sequence of mass ratios that support solitary waves in the $1\!:\!1$ \cite{jvs-solitary0} and $1\!:\!2$ \cite{jvs-solitary1} dimers with Hertzian spring forces.
With Gendelman, they find in the $1\!:\!1$ dimer a different sequence of mass ratios tending to 0 for which waves asymptote to oscillatory pulses \cite{jsvg}.

\subsection{Numerical method}
The main technical problem that we face in this paper is that standard path-continuation software packages such as AUTO \cite{doedel1998auto} and PDE2Path \cite{uecker2012pde2path} cannot be applied to \eqref{tw eqns} on account of the shifts in the arguments. 
Early numerical work involving MFDEs can be found in \cite{VL7}, which was continued by Elmer and Van Vleck in the extensive series of papers \cite{EVVRCW,EVVSDRD,EVVACTW,EVVPerDiff,EVV}. 
Results specific to FPUT-type problems can be found in \cite[{\S}4]{faver-goodman-wright}, where the authors use Fourier decompositions to attack the shifted terms.

Our computations here involve the use of a collocation solver based on \cite{EVVRCW,HJHVL2005} that is able to solve MFDEs on finite intervals. 
In particular, it can handle general $n$-component problems of the form
\begin{equation}
\label{int:eq:mfdeCode}
\phi'(\xi) = f\Big( \phi(\xi), \phi\big(\xi + \sigma_1(\xi, \phi(\xi)) \big), \ldots , \phi\big(\xi + \sigma_N(\xi, \phi(\xi))\big) \Big),
\end{equation}
for given functions $f: \R^{n(N+1)} \to \R^n$, $\tau: \R \to \R^n$ and shifts $\sigma_i: \R^{1+n} \to \R$.
This is achieved by representing $\phi(\xi)$ on each grid-interval in terms of a standard Runge-Kutta monomial basis, requiring \eqref{int:eq:mfdeCode} to be satisfied at each of the interior Gaussian collocation points.

Various types of boundary conditions can be used to close the system \eqref{int:eq:mfdeCode}, which we exploit heavily here in order to ensure that our computed waves have the required ``solitary + ripple'' structure.
Notice that  the shifts $\sigma_i$ may depend on the spatial variable $\xi$ as well as on the function value $\phi(\xi)$ itself.
This allows us to compute periodic solutions to MFDEs even when the period is unknown, which is essential for our purposes here.

\subsection{Outline}
In {\S}\ref{sec:bck} we introduce our computational coordinate system that respects certain important symmetries. 
In addition, we summarize the mathematical background behind Beale's decomposition procedure for our diatomic waves. 
We numerically analyze two scalar MFDEs in {\S}\ref{sec:mono} that are related to the monatomic limit, which allows us to compute the leading-order ripple amplitude of our micropterons. 
The solutions to the full diatomic wave MFDE \eqref{tw eqns} are computed in {\S}\ref{sec:di}, while {\S}\ref{sec:full:fpu} describes our direct simulations of the original FPUT problem \eqref{newton}.
We close in {\S}\ref{sec:disc} with a brief discussion of possible future research directions.

\section{Background}\label{sec:bck}
Our goal here is to briefly outline the mathematical background required to appreciate the choices made during our numerical work in later sections. 
In order to streamline our presentation, we start by introducing the exponentially localized and periodic Sobolev spaces
\begin{equation}
H^r_q 
:= \set{f \in H^r}{\cosh^q(\cdot)f \in H^r},
\qquad
H^r_{\per} 
:= \set{f \in H^r([0,2\pi])}{f(0) = f(2\pi)},
\end{equation}
using the natural norm
\begin{equation}
\norm{f}_{H_q^r}
:= \norm{\cosh^q(\cdot)f}_{H^r}
\end{equation}
for the weighted spaces.
In addition, we define the odd subspaces
\begin{equation}
O^r_q 
:= H^r_q \cap \{\text{odd functions}\},
\quad
O^r_{\per} 
:= H^r_{\per} \cap \{\text{odd functions}\},
\end{equation}
together with their even counterparts
\begin{equation}
E^r_q 
:= H^r_q \cap \{\text{even functions}\},
\qquad
E^r_{\per}
:= H^r_{\per} \cap \{\text{even functions}\}.
\end{equation}
At times, we restrict the latter even further and consider the ``mean-zero'' spaces
\begin{equation}
E^r_{q,0} 
:= \set{f \in E^r_q}{\int_{-\infty}^\infty f(\xi) \dxi = 0},
\qquad 
E^r_{\per,0} 
:= \set{f \in E^r_{\per}}{\int_0^{2\pi} f(\xi) \dxi = 0}.
\end{equation}

\subsection{Coordinate system}
In order to exploit several useful symmetries, it is convenient to introduce new variables for the traveling wave problem. 
In particular, we introduce a new mass parameter $\mu = 1/m - 1$ to encode the deviation from the equal-mass limit and consider the linear combinations
\begin{equation}\label{sbar-rbar cov}
\sbar_1 =\frac{\rbar_o + \rbar_e}{2}
\quadword{and}
\sbar_2 = \frac{\rbar_o - \rbar_e}{2}.
\end{equation}
Writing $\sbar = (\sbar_1,\sbar_2)$,
the diatomic traveling wave MFDE \eqref{tw eqns} can be recast into the form
\begin{equation}\label{eq:trv;wave:s1:s2}
-c^2 \sbar'' = \D_{1}(\mu)\begin{pmatrix*}
\sbar_1 + \sbar_1^2 + \sbar_2^2 \\
\sbar_2 + 2  \sbar_1 \sbar_2
\end{pmatrix*},
\end{equation}
Here we have introduced the linear operator
\begin{equation}
\D_{\nu}(\mu)
:= \frac{1}{2}
\begin{bmatrix*}
(2 + \mu)(2 - A_{\nu}) & \mu \delta_{\nu} \\
- \mu \delta_{\nu} & (2 + \mu)(2 + A_{\nu})
\end{bmatrix*},
\end{equation}
in which $A_{\nu}$ and $\delta_{\nu}$ are the shift operators
\begin{equation}\label{A delta defns}
[A_{\nu}f](\xi) := f(\xi + \nu) + f(\xi - \nu)
\quadword{and}
[\delta_{\nu}f](\xi) := f(\xi + \nu) - f(\xi - \nu).
\end{equation}

It is easy to check that the system \eqref{eq:trv;wave:s1:s2} preserves an ``even $\times$ odd'' symmetry. 
More precisely,
if $\sbar_1$ is even and $\sbar_2$ is odd, then the first components of both sides of \eqref{eq:trv;wave:s1:s2} are even (and also mean-zero), while the second components are both odd.
A second important symmetry can be readily observed in our original traveling wave problem \eqref{tw eqns}. 
Indeed, this system is invariant under the transformation
\begin{equation}
\label{eq:trnsf:m:vs:m:inv}
(\rbar_o, \rbar_e, m, c)
\mapsto (\rbar_e, \rbar_o, 1/m, c\sqrt{m}).
\end{equation}
We exploit this to perform our numerics on the symmetrized system \eqref{eq:trv;wave:s1:s2} in the bounded regime $\mu \in (-1,0]$, which corresponds to $m \in [1, \infty)$. 
Using \eqref{eq:trnsf:m:vs:m:inv} we subsequently transfer these computations back to the regime $m \in (0,1]$, which we feel is much better suited for discussing and visualizing our results.

\subsection{Periodic waves}
Here we summarize the procedure used in \cite{faver-wright, hoffman-wright, faver-hupkes} to construct the periodic traveling wave solutions that constitute the background ripples of our diatomic waves.
We make a few minor changes to the parametrization used in these papers that will facilitate various parts of our subsequent numerics.

\subsubsection{The linearized periodic problem}\label{lin per section}
We start by looking for solutions to the linearization of \eqref{eq:trv;wave:s1:s2}, which is 
\begin{equation}\label{symm tw lin}
-c^2\sbar'' 
= \D_1(\mu)\sbar.
\end{equation}
We seek solutions of the form $\sbar(\xi) = e^{i \omega \xi}\sbar_k$, where $\sbar_k \in \C^2$.
We find that the vector $\sbar_k$ must satisfy the characteristic relation
\[
\Delta(\omega;c,\mu)\sbar_k
= 0,
\]
where
\begin{equation}
\Delta(\omega;c,\mu) 
:= \begin{bmatrix*}
-c^2 \omega^2 + (2 + \mu)\big(1 - \cos(\omega)\big) & i \mu \sin(\omega) \\
- i \mu \sin(\omega) & -c^2 \omega^2 +  (2 + \mu)\big(1 + \cos(\omega)\big)
\end{bmatrix*}.
\end{equation}
Upon introducing the expressions
\begin{equation}
\lambda_{\mu}^{\pm}(\omega) 
:= 2 + \mu \pm \sqrt{\mu^2 + 4 (1 + \mu) \cos^2(\omega) },
\end{equation}
together with
\begin{equation}
\mathcal{B}_{\pm}(\omega;c,\mu) 
:= - c^2 \omega^2 + \lambda_{\mu}^{\pm}(\omega),
\end{equation}
we have the convenient factorization
\begin{equation}
\det\big(\Delta(\omega;c,\mu)\big) 
= \mathcal{B}_-(\omega;c,\mu) \mathcal{B}_+(\omega;c,\mu).
\end{equation}

For $\mu \ge -1$  we have the useful inequality\footnote{
Note that $C_{\mu}$ corresponds with the critical
speed $c^{(s)}_{m}$ defined in \eqref{lw speed of sound}.
} 
\cite[Eq.\@ (C.1.5)]{faver-hupkes}
\begin{equation}
|(\lambda_{\mu}^{\pm})'(\omega)| 
\le 2C_{\mu}^2|\omega|,
\qquad C_{\mu} := \sqrt{\frac{2(1+\mu)}{2+\mu}}.
\end{equation}
For $c$ and $\mu$ satisfying $c > C_{\mu}$, this allows us to conclude that $\mathcal{B}_-$ has no zeros other than $\omega = 0$, with the corresponding eigenvector $\Delta(0;c,\mu)(1,0)^{\transpose} =0$. 
In addition, the observations
\begin{equation}
\mathcal{B}_+(0;c,\mu) > 0,
\qquad
\mathcal{B}_+(\infty;c,\mu) = -\infty,
\qquad
\mathcal{B}_+'\big((0, \infty) ;c,\mu\big) < 0 
\end{equation}
imply that there is a unique $\omega_{c,\mu} > 0$ with $\mathcal{B}_+(\omega_{c,\mu};c,\mu) = 0$.
The corresponding eigenvector
\begin{equation}
\Delta(\omega_{c,\mu};c,\mu)
\big(\nu_1^{c,\mu}, \nu_2^{c,\mu}\big)^{\transpose}
= 0
\end{equation}
can be chosen to be continuous in $(c,\mu)$, while satisfying the normalization
\begin{equation}
\big|\nu_1^{c,\mu}\big|^2 + \big|\nu_2^{c,\mu}\big|^2 = 1
\quadword{and}
\big(\nu_1^{c,0}, \nu_2^{c,0}\big) = (0,1).
\end{equation}
The corresponding solution to the linearization \eqref{symm tw lin} is then given by
\begin{equation}
\sbar_{\lin}^{c,\mu}(\xi) 
:= \big(\nu_1^{c,\mu} \cos(\omega_{c,\mu} \xi), \nu_2^{c,\mu} \sin(\omega_{c,\mu} \xi)\big).
\end{equation}

\subsubsection{The full nonlinear periodic problem}\label{full nl periodic section}
One can subsequently construct solutions for the full nonlinear periodic problem \eqref{eq:trv;wave:s1:s2} via a Crandall-Rabinowitz-Zeidler ``bifurcation from a simple eigenvalue'' argument \cite{crandall-rabinowitz, zeidler}.
This provides triplets
\begin{equation}
\label{eq:bck:triplet:p:k}
\big(\pbar_1^{c,\mu}[a], \pbar_2^{c,\mu}[a], \omega_{c,\mu}[a]\big) \in E_{\per,0}^2 \times O_{\per}^2 \times \mathbb{R}
\end{equation}
that are parametrized by the (small) signed amplitude
\begin{equation}
\label{eq:per:normalization}
a 
= \mathrm{sign}\big(\pbar_1^{c,\mu}[a](0)\nu_1^{c,\mu} + (\pbar_2^{c,\mu}[a])'(0)\nu_2^{c,\mu}/\omega_{c,\mu}\big)\sqrt{\norm{\pbar_1^{c,\mu}[a]}_{L^{\infty}}^2 + \norm{\pbar_2^{c,\mu}[a]}_{L^{\infty}}^2}
\end{equation}
and that yield solutions to 
\eqref{eq:trv;wave:s1:s2} of the form
\begin{equation}
\sbar(\xi) 
= \pbar^{c,\mu}[a](\omega_{c,\mu}[a] \xi)
:= \big(\pbar_1^{c,\mu}[a],\pbar_2^{c,\mu}[a]\big) (\omega_{c,\mu}[a] \xi).
\end{equation}
We note that the zero eigenvalue is ruled out by taking $\pbar_1^{c,\mu}[a]$ to be a mean-zero function. 
See \cite{faver-wright, hoffman-wright, faver-hupkes} for the details of this construction in the various parameter regimes.

These solutions branch off from the linearized solution in the sense that
\begin{equation}
\pbar^{c,\mu}[a] = a\sbar_{\lin}^{c,\mu} + \O(a^2)
\quadword{and}
\omega_{c,\mu}[a] = \omega_{c,\mu} + \O(a).
\end{equation}
Indeed, we view the $\mathrm{sign}(\cdot)$ factor in \eqref{eq:per:normalization} as defining an orientation on $a$ relative to $\sbar_{\lin}^{c,\mu}$. 
Notice in addition that $|a| = \norm{\pbar^{c,\mu}[a]}_{L^{\infty}}$, which implies that $|a|$ is genuinely the ``amplitude'' of the periodic profile $\sbar$, which is well-suited for our purposes here. 

We emphasize that our explicit choice \eqref{eq:per:normalization} differs slightly from the definitions in \cite{faver-wright, hoffman-wright, faver-hupkes}, where the $L^{\infty}$-norm of the periodic profiles is only $\O(a)$ up to some $c$-dependent constants. 
The equivalence of these different approaches to $a$ follows from the fact that we have Lipschitz-smoothness in $a$, see Lemma C.1 in \cite{hoffman-wright} or part (ii) of Lemma C.3 in \cite{faver-hupkes}.
We discuss our numerical simulations of these periodic solutions in {\S}\ref{sec:periodic numerics}.

\subsection{Micropterons}
Here we fix $|c| \gtrsim 1$ and consider masses $m \approx 1$, which allows us to use $\mu$ as a small parameter. We set out to find solutions to \eqref{eq:trv;wave:s1:s2} in the vicinity of the pair
\begin{equation}
\label{eq:bck:microp:pair:s1:s2}
(\sbar_1, \sbar_2) 
= (\phi_c, 0),
\end{equation}
where $\phi_c$ solves the monatomic traveling wave problem \eqref{mono:fput:tw:eqn}.

\subsubsection{The linearization at $(\phi_c,0)$}
\label{sec:sub:lin:phi:c}
The linearization of our symmetrized traveling wave problem \eqref{eq:trv;wave:s1:s2} around the pair \eqref{eq:bck:microp:pair:s1:s2} at $\mu = 0$ is the diagonal operator $\diag(\H_c,\L_c)$, whose components are given by
\begin{equation}\label{Hc Lc defn}
\H_c v_1 := c^2 v_1'' + (2 - A_1) \big[ v_1 + 2 \phi_c v_1 \big]
\quadword{and}
\L_c v_2 = c^2 v_2''+ (2 + A_1) \big[ v_2 + 2 \phi_c v_2 \big].
\end{equation}

The shift operator $A_1$ was defined in \eqref{A delta defns}.
The first component $\H_c$ is the linearization of the monatomic traveling wave problem \eqref{mono:fput:tw:eqn} at $\phi_c$.
For $c \gtrsim 1$ it is known \cite{hoffman-wright, faver-hupkes} that this operator is invertible  from $E^{r+2}_q$ to $E_{q,0}^r$ for suitably small $q$ and all $r \ge 0$.

The second component $\L_c$ was analyzed in detail in \cite{faver-hupkes}, where it is shown that $\L_c$ is injective from $O^{r+2}_q$ into $O^r_q$ with
\begin{equation}\label{eq:char:range:L}
\mathrm{Range}(\L_c) 
= \set{g \in O^r_q}{\ip{g}{\gamma_c}_{L^2} = 0}.
\end{equation}
Here the odd bounded function $\gamma_c$ is a nontrivial solution to the adjoint problem
\begin{equation}\label{eq:bck:mfde:for:jost:gamma}
-c^2 \gamma_c'' 
= (1 + 2 \phi_c) (2 + A_1) \gamma_c.
\end{equation}
Recalling the critical frequency $\omega_{c,0}$ from {\S}\ref{lin per section}, it has the limiting behavior
\begin{equation}
\lim_{\xi \to \infty} |\gamma_c(\xi) - \sin\big(\omega_{c,0}(\xi + \vartheta_c) \big)|
= 0
\end{equation}
for some asymptotic phase-shift $\vartheta_c$.
More specifically, we have the decomposition
\begin{equation}
\gamma_c(\xi) = \upsilon_c(\xi) + \sin\big(\omega_{c,0}(\xi + \vartheta_c) \big)
\end{equation}
in which $\upsilon_c$ is exponentially localized.
We call such a function $\gamma_c$ a Jost solution for the adjoint problem \eqref{eq:bck:mfde:for:jost:gamma}, in the spirit of classical Jost solutions for the Schr\"{o}dinger operator.
We outline our numerical procedure for the computation of these Jost solutions in {\S}\ref{sec:jost numerics}.

The asymptotic frequency $\omega_{c,0}$ and the asymptotic phase shift $\vartheta_c$ interact in a special way.
Recalling the branch \eqref{fr-p speed of sound}--\eqref{fp lw limit} of solitary waves constructed by Friesecke and Pego, a careful bifurcation analysis that combines MFDE theory with residue calculus \cite[Eq.\@ (6.3.14), (6.4.11), App.\@ E.3.2]{faver-hupkes} delivers the expansion
\begin{equation}
\label{eq:exp:phase:shift}
\omega_{\cep,0} \vartheta_{\cep}
=  - \frac{ \ep \omega_{1,0}^2}{8\mathcal{B}_+'(\omega_{1, 0}; 1, 0)} \int_{-\infty}^{\infty} \mathrm{sech}^2(\xi/2) \dxi + \O(\ep^2).
\end{equation}
This implies that $\sin(\omega_{\cep,0} \vartheta_{\cep}) \neq 0$ for small $\ep > 0$, which is essential in the discussion below.

\subsubsection{Beale's ansatz}
Following Beale \cite{beale}, we now search for solutions to the traveling wave problem \eqref{eq:trv;wave:s1:s2} that have the form
\begin{equation}
\label{eq:bck:beale:ansatz}
(\sbar_1,\sbar_2)(\xi)
= (\phi_c, 0)(\xi)
+ (w_1, w_2)(\xi) 
+ \big(\pbar_1^{c,\mu}[a],\pbar_2^{c,\mu}[a]\big)(\xi) ,
\end{equation}
where we take
\begin{equation}
(w_1, w_2, a) \in E^2_q \times O^2_q \times \R
\end{equation}
with $q > 0$ suitably small.
One can use the invertibility of $\H_c$ from \eqref{Hc Lc defn} to construct a fixed point equation for the unknown ``error'' $w_1$, but obtaining equations for $w_2$ and the ``amplitude'' $a$ is more challenging.

The system's second component can be informally written as
\begin{equation}
\L_c w_2
+ \mu\chi_c + a\eta_c 
=  \O \big( a^2 + \mu^2 + \norm{w_1}_{H^2_q}^2 +\norm{w_2}_{H^2_q}^2\big).
\end{equation}
Since the periodic profiles $(\pbar_1^{c,\mu}[a], \pbar_2^{c,\mu}[a])$ already solve \eqref{eq:trv;wave:s1:s2}, the function $\eta_c$ can be recognized as the cross-term
\begin{equation}
\eta_c
= (2+A_1)\big[\phi_c \sin(\omega_{c,0}\cdot)\big].    
\end{equation}
In addition, since $(\phi_c, 0)$ satisfies \eqref{eq:trv;wave:s1:s2} at $\mu=0$, the function $\chi_c$ arises by applying the difference $\D_1(\mu) -\D_1(0)$ to this pair, yielding
\begin{equation}
\chi_c
= - \frac{1}{2} \delta_1\big[ \phi_c + \phi_c^2\big].
\end{equation}

\subsubsection{The ripple amplitude revealed}\label{sec: micropt ampl}
The range characterization \eqref{eq:char:range:L} now readily yields
\begin{equation}
\label{eq:bck:expr:coeff:ampl}
a = - K_c \mu  + \O(\mu^2),
\qquad  \hbox{ with } \qquad K_c := \frac{\ip{\gamma_c}{\chi_c}_{L^2}}{\ip{\gamma_c}{\eta_c}_{L^2}},
\end{equation}
which is well-defined for $c = \cep$ and small $\ep > 0$ on account of the remarkable explicit (corrected) identity
\cite[Eq. 5.3.7]{faver-hupkes}
\begin{equation}
\label{eq:bck:expr:int:chi:c:a}
\ip{\gamma_c}{\eta_c}_{L^2}
= - \mathcal{B}_+'(\omega_{c,0};c,0)\sin(\omega_{c,0}\vartheta_c)    
\end{equation}
and the remark following the expansion \eqref{eq:exp:phase:shift}.

If $K_c \ne 0$, then the leading-order coefficient for $a$ in \eqref{eq:bck:expr:int:chi:c:a} is nonzero, implying that our ripple's amplitude is only algebraically small in $\mu$.
That is, we will have produced a genuine micropteron.
It is far from clear, however, whether or not $K_c$ can ever vanish for a certain choice of $c$.
Taking $c=\cep$ and $\phi_{\cep}$ to be a near-sonic Friesecke-Pego solitary wave, we might attempt an expansion of $K_{\cep}$ in powers of $\ep$, since there are a host of tight $\ep$-estimates on $\phi_{\cep}$ \cite{friesecke-pego1,hoffman-wayne, faver-hupkes}.
However, we conjecture that such an expansion will actually reveal that $K_{\cep}$ is small beyond all algebraic orders of $\ep$.
Consequently, one of the computational goals of this paper is to numerically investigate $K_c$; see {\S}\ref{sec:mono}.

\subsection{Nanopterons}
In the long wave \cite{faver-wright} and small mass \cite{hoffman-wright} limits,
the traveling waves are nanopterons, and so the amplitude $a$ of the ripples is small beyond all orders of the relevant small parameter.
There is no question, then, of attempting to isolate its leading order behavior as we do for the micropteron.
In particular, we do not attempt to compute the analogous Jost solutions that appear in those problems as solutions to certain auxiliary MFDEs.

We do discuss, however, how the formal $m=0$ solitary waves $\varphi_c$ from {\S}\ref{small mass section intro} behave in the $\sbar_1$ and $\sbar_2$ coordinates.
Let $\varphi_c$ be an even solitary wave solution of the MFDE \eqref{small mass MFDE}.
Since this MFDE is shift-invariant, the profile $\tilde{\varphi}_c(\xi) := \varphi_c(\xi+1/4)$ is also a solution.
The identities \eqref{eq:id:small:mass} and \eqref{small mass big R} tell us that putting
\[
\rbar_o(\xi)
= \tilde{\varphi}_c(\xi/2)
\quadword{and}
\rbar_e(\xi)
= \tilde{\varphi}_c\big((\xi-1)/2\big)
\]
formally solves the original FPUT traveling wave problem \eqref{tw eqns} at $m=0$.
Then using the change of variables \eqref{sbar-rbar cov}, we find that these solutions read
\[
\sbar_1(\xi)
= \frac{\varphi_c\big((\xi +1/2)/2\big)+\varphi_c\big((\xi -1/2)/2\big)}{2}
\]
and
\[
\sbar_2(\xi)
= \frac{\varphi_c\big((\xi +1/2)/2\big)- \varphi_c\big((\xi -1/2)/2\big)}{2}.
\]
In particular - in contrast to the equal-mass limit - the $\sbar_2$ component does \textit{not} vanish in the limit $m \downarrow 0$; see, for example, graph III in Fig.\@ \ref{fig:cross:c:1p65}.

\section{Monatomic simulations}\label{sec:mono}
Our goal here is to numerically find solutions to the monatomic traveling wave MFDE \eqref{mono:fput:tw:eqn}, together with the problem \eqref{eq:bck:mfde:for:jost:gamma} that describes the associated Jost solutions.
Compared to the discussion in {\S}\ref{sec:bck} above, the main change here is that we do not use the wave speed\footnote{
From now on we use $\sigma$ for the wave speed in order to emphasize that it is a numerically computed variable
and not a fixed system parameter.
} 
$\sigma > 1$ to parametrize the waves, but rather introduce a new parameter
\begin{equation}
\label{eq:mono:def:kappa}
\kappa \sim \sqrt{8\phi(0)}
\end{equation}
that is a measure for the center amplitude of the wave. 
This parameter is always supplied {\it{a-priori}} to our numerical method.

The factor of 8 arises from the estimates for the Frisecke-Pego monatomic solutions in \eqref{fp lw limit}. 
Indeed, in view of \cite[Eq.\@ (4.1)]{friesecke-pego1} this can be seen as an equivalent parametrization for the family \eqref{fr-p speed of sound}--\eqref{fp lw limit}, with the convenient property
\begin{equation}
\label{eq:mono:eqv:ep:kappa}
\ep = \kappa + \O(\kappa^2).
\end{equation}
The numerical motivation for this choice is that it is much easier to use fixed boundary conditions for the wave profile. 
This helps to stabilize the code and also --- crucially --- prevents convergence to the omnipresent zero solution. 

\begin{figure}[t]
\begin{center}
\includegraphics[width=1\textwidth]{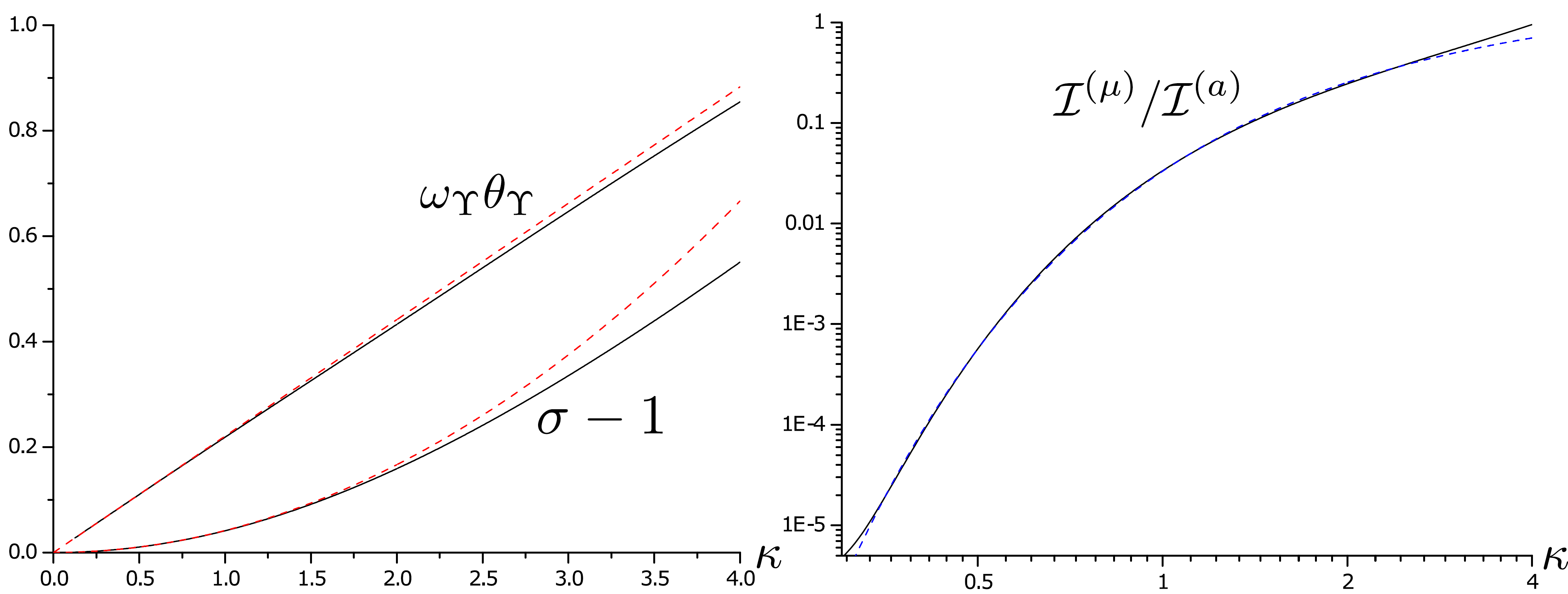}
\caption{The left graph visualizes the computed values for the asymptotic phase-shift of the Jost solution and the monatomic wave speed (solid black), together with their leading-order predictions \eqref{eq:mono:predictions} (dashed red). The right graph contains our computed values for $-K_{\sigma}$
(solid black),
together with the fit \eqref{eq:mono:pred:K:c} (dashed blue).
} 
\label{fig:mono:coeffs}
\end{center}
\end{figure}

\subsection{Wave speed and profile}
For fixed $\kappa > 0$, we aim to find a solution
\begin{equation}
\phi(\xi) \sim \kappa^2 \Phi(\kappa \xi),
\end{equation}
to \eqref{mono:fput:tw:eqn} with $c = \sigma$ by numerically computing the pair $(\sigma, \Phi)$. 
This scaling is inspired by the limiting behavior \eqref{fp lw limit} and allows us to use the same numerical interval for a wide range of values of $\kappa$. 
This enables us to use a continuation approach where we gradually modify $\kappa$, without the danger that our solutions become too wide. 

In order to find our numerical wave, we introduce a computational coordinate
\begin{equation}
\tau = \kappa \xi \in [0, L]
\end{equation}
for some fixed $L > 0$ and set out to solve the problem
\begin{equation}
\label{eq:mono:fde:phi}
- \kappa^2 \sigma^2 \Phi''(\tau) = \Big[(2 - A_{\kappa}) \big( \Phi + \kappa^2 \Phi^2 \big)\Big](\tau),
\end{equation}
augmented by the boundary conditions
\begin{equation}
\label{eq:mono:bnd:phi}
\Phi(0) = \frac{1}{8},
\qquad \qquad
\Phi'(0) = 0,
\qquad \qquad
\Phi(L) = 0.
\end{equation}
The first condition is related to \eqref{eq:mono:def:kappa}, while the second allows us to (virtually) extend $\Phi$ to an even function.
In particular, we take $\Phi(\tau) = \Phi(-\tau)$ whenever the shifts in \eqref{eq:mono:fde:phi} require an evaluation of $\Phi$ at a negative argument.

\subsection{Jost solutions}\label{sec:jost numerics}
In order to solve \eqref{eq:bck:mfde:for:jost:gamma}, we write
\begin{equation}
\label{eq:mono:ansatz:gamma:sigma}
\gamma_{\sigma}(\xi) 
\sim
\sin\big(  \omega_\Upsilon (\xi + \theta_\Upsilon) \big) 
+
\beta_\Upsilon \Upsilon(\kappa \xi)
\end{equation}
and set out to numerically compute $(\omega_{\Upsilon}, \theta_{\Upsilon}, \beta_{\Upsilon}, \Upsilon)$. 
The first of these is the solution to the scalar nonlinear problem
\begin{equation}
\sigma^2 \omega_{\Upsilon}^2 = 2 +2 \cos(\omega_{\Upsilon}).
\end{equation}
The remainder function $\Upsilon$ is characterized by the MFDE
\begin{equation}
\begin{array}{lcl}
-\kappa^2 \sigma^2  \Upsilon''(\tau) & = & 
\Big( 1+ 2 \kappa^2 \Phi(\tau)\Big)
\Big( 2 \Upsilon(\tau) + [A_{\kappa}\Upsilon](\tau)\Big)
\\[0.2cm]
& & \qquad
+ 2 \kappa^2  \sigma^2 \omega_\Upsilon^2 \beta_\Upsilon^{-1} \Phi(\tau) \sin\big( \omega_\Upsilon (\tau/\kappa + \theta_\Upsilon) \big),
\end{array}
\end{equation}
augmented by the boundary conditions
\begin{equation}
\label{eq:mono:bnd:upsilon}
\Upsilon(0) = - \beta_\Upsilon^{-1} \sin(\omega_\Upsilon \theta_\Upsilon), \qquad \Upsilon'(0) =1,
\qquad \Upsilon(L) = \Upsilon'(L) = 0.
\end{equation}
The first of these ensures that 
$\gamma_{\sigma}$ can be (virtually) extended to an odd function by writing
\begin{equation}
- \Upsilon(-\tau)
= \Upsilon(\tau) + \beta_\Upsilon^{-1} \sin\big(\omega_\Upsilon (\tau/\kappa+\theta_\Upsilon) \big)
+  \beta_\Upsilon^{-1} \sin\big(\omega_\Upsilon( - \tau/\kappa +  \theta_\Upsilon) \big).
\end{equation}
The second boundary condition is a convenient normalization,
but it does require us to introduce the extra parameter $\beta_\Upsilon$ in the ansatz \eqref{eq:mono:ansatz:gamma:sigma}. 

\subsection{The $K_c$ coefficient}
In order to evaluate the explicit coefficient in \eqref{eq:bck:expr:coeff:ampl}, we introduce the function
\begin{equation}
\Psi^{(\eta)}(\tau)
=  \Phi(\tau+\kappa) \sin( \omega_{\Upsilon} \tau/\kappa + 1)
+ 2 \Phi(\tau)
\sin(\omega_{\Upsilon}\tau/\kappa)
+ \Phi(\tau-\kappa) \sin( \omega_{\Upsilon} \tau/\kappa - 1),
\end{equation}
together with
\begin{equation}
\Psi^{(\chi)}(\tau)
= - \frac{1}{2}
\Big[\Phi(\tau + \kappa) + \kappa^2 \Phi(\tau + \kappa)^2
- \Phi(\tau - \kappa)
-\kappa^2 \Phi(\tau - \kappa)^2
\Big].
\end{equation}
This allows us to define the integrals
\begin{equation}
\label{eq:mono:itgs}
\begin{array}{lcl}
\mathcal{I}^{(\eta)} & = & 
2 \kappa \int_0^L 
\big[ \sin\big(\omega_{\Upsilon}(\tau/\kappa+\theta_{\Upsilon})\big)
+ \beta_\Upsilon \Upsilon(\tau)
\big] \Psi^{(\eta)}(\tau)
\, d \tau ,
\\[0.2cm]
\mathcal{I}^{(\chi)} &=& 
2 \kappa \int_0^L \big[ \sin\big(\omega_{\Upsilon}(\tau/\kappa+\theta_{\Upsilon})\big)
+ \beta_\Upsilon \Upsilon(\tau)
\big] \Psi^{(\chi)}(\tau)
\, d \tau ,
\end{array}
\end{equation}
which should be seen as our numerical proxies for the inner products $\langle \gamma_c, \eta_c \rangle_{L^2}$ respectively $\langle \gamma_c, \chi_c \rangle_{L^2}$ that appear in \eqref{eq:bck:expr:coeff:ampl}. 
In particular, we obtain the numerical prediction
\begin{equation}
\label{eq:mono:ampl:coeff}
K_{\sigma} \sim - \mathcal{I}^{(\chi)}/ \mathcal{I}^{(\eta)}.
\end{equation}
We note that the identity \eqref{eq:bck:expr:int:chi:c:a} implies that we expect to have
\begin{equation}
\mathcal{I}^{(\eta)} \approx - \mathcal{B}_+'(\omega_{\Upsilon};\sigma,0)\sin(\omega_{\Upsilon} \theta_{\Upsilon}),
\end{equation}
which we used as an independent monitor for the accuracy of our numerical schemes. 

\subsection{Implementation}
Note that the two scalar differential equations  for $\Phi$ and $\Upsilon$ are both of order two, while there are three free parameters $(\sigma, \beta_\Upsilon, \theta_\Upsilon)$ that need to be determined. 
The collocation solver discussed in \cite{HJHVL2005} hence requires seven boundary conditions, which indeed matches the number supplied in \eqref{eq:mono:bnd:phi} and \eqref{eq:mono:bnd:upsilon}.
We solved the combined system on the interval $[0,L] = [0,32]$ for a range of values for $\kappa \ge 1/8$.
The results for $\sigma$ and the product $\omega_{\Upsilon} \theta_{\Upsilon}$ can be found in Fig.\@ \ref{fig:mono:coeffs}a.

The integrals \eqref{eq:mono:itgs} were computed by applying the mid-point rule with $10^6$ gridpoints. 
The resulting values can be found in Fig.\@ \ref{fig:mono:coeffs}b.
Due to the high-frequency oscillations in the integrand for $\mathcal{I}^{(\chi)}$ that appear as $\kappa \downarrow 0$, the values for this integral become unreliable when $\kappa$ is too small. 
For this reason, we restricted the plot to the range $\kappa \ge 0.3$. 
This cut-off was determined by changing the number of gridpoints used for the integral evaluation and checking whether the computed values remain stable.

\subsection{Discussion}
Evaluating \eqref{fr-p speed of sound}, and \eqref{eq:exp:phase:shift} using the reparametrization \eqref{eq:mono:eqv:ep:kappa} and the observation $\omega_{1,0} \sim 1.478170266$, we arrive at the predictions
\begin{equation}
\label{eq:mono:predictions}
\omega_\Upsilon \vartheta_\Upsilon
\sim 0.2208053960 \kappa + \O(\kappa^2),
\qquad \qquad
\sigma \sim 1 + \frac{1}{24} \kappa^2 + O (\kappa^3).
\end{equation}
These predictions for the phase-shift and wave speed agree remarkably well with our numerics; see Fig.\@ \ref{fig:mono:coeffs}a.
We also emphasize that we are able to find waves for relatively large values of $\kappa$, which (arguably) fall outside of the small-amplitude regime analyzed in \cite{friesecke-pego1,friesecke-wattis}.

As discussed in {\S}\ref{sec: micropt ampl}, we have no {\it{a-priori}} predictions for the amplitude coefficient \eqref{eq:mono:ampl:coeff}.
To examine this in further detail, we fitted the graph with a function that is exponentially small in $\kappa$ and found the approximation
\begin{equation}
\label{eq:mono:pred:K:c}
\mathcal{I}^{(\chi)}/ \mathcal{I}^{(\eta)}
\approx
1.93756 \, \mathrm{exp}\big[-4.06704 / \kappa \big],
\end{equation}
which agrees quite well with our data points.
In particular, we view this as support for our conjecture that $K_{c_{\epsilon}}$ is small beyond all orders in $\epsilon$. 
In any case, we now have solid evidence to show that $K_{\sigma} < 0$ for a wide range of $\sigma > 1$, which in our view justifies the micropteron terminology that we used in \cite{faver-hupkes}.

\section{Diatomic Simulations}
\label{sec:di}

\begin{figure}[t]
\begin{center}
\includegraphics[width=1\textwidth]{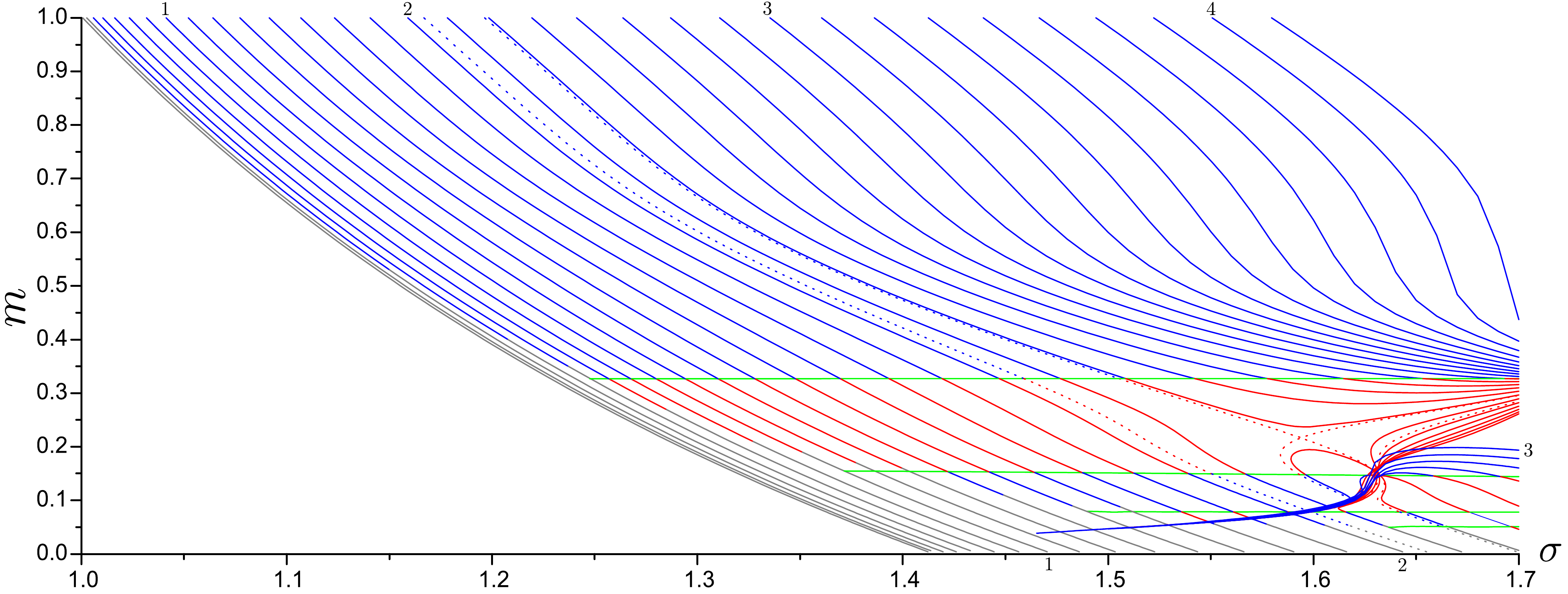}
\caption{Overview of several iso-$\kappa$ curves where diatomic waves exist. The small integers denote the $\kappa$ values, while the dashed separatrices correspond to the special values \eqref{eq:di:spec:kappa}. The color codes are described in {\S}\ref{sec:di:impl}. Notice that the green branches of solitary waves and the blue fold curve terminate for technical reasons, but in principle extend further into the grey ``small-ripple'' regime.
} 
\label{fig:di:overview}
\end{center}
\end{figure}

We are now ready to search for solutions of the form \eqref{eq:bck:beale:ansatz} to the full diatomic wave problem \eqref{eq:trv;wave:s1:s2}.
We reuse the scaling parameter $\kappa$, which now should be interpreted as
\begin{equation}\label{kappa approx sol core}
\kappa \sim \sqrt{8\big(\phi_{\sigma}(0) + w_1(0) \big) }. 
\end{equation}
In particular, this parameter is a measure for the center amplitude of the \textit{solitary} part of the wave, which includes both the known solitary core in Beale's ansatz \eqref{eq:bck:beale:ansatz} as well as the unknown localized terms, but excludes the background periodic ripple.
Since the size of this ripple is zero at $m = 1$ and $m=0$, this coincides with the parameter $\kappa$ that we used in the monatomic setting of {\S}\ref{sec:mono}.

\subsection{Periodic solutions}\label{sec:periodic numerics}

We first aim to construct candidates for the background ripple by searching for periodic solutions to \eqref{eq:trv;wave:s1:s2} of the form
\begin{equation}
(\sbar_1, \sbar_2)(\xi) 
\sim 
\beta_P (\tilde{P}_1, \tilde{P}_2)( \omega_P \xi ).
\end{equation}
In particular, in terms of the computational coordinate
\begin{equation}
\tilde{\tau} = \omega_P \xi \in [0,L],
\end{equation}
we need to solve the problem 
\begin{equation}
\begin{array}{lcl}
-\omega_P^2 \sigma^2 \tilde{P}'' =
\D_{\omega_P}(\mu)
\begin{pmatrix*}
\tilde{P}_1 + \beta_P (\tilde{P}_1^2 + \tilde{P}_2^2) \\[0.2cm]
\tilde{P}_2 + 2 \beta_P \tilde{P}_1 \tilde{P}_2
\end{pmatrix*},
\end{array}
\end{equation}
augmented by the (boundary) conditions
\begin{equation}
\label{eq:di:per:bnd}
\int_{0}^L \tilde{P}_1(\tilde{\tau}) \, d \tilde{\tau} = 0,
\qquad
\tilde{P}_1'(0) = 0,
\qquad
\tilde{P}_2(0) = \tilde{P}_2(L) = 0
\end{equation}
and the normalization 
\begin{equation}
\tilde{P}_1(0)^2 +  \tilde{P}_2'(0)^2 = 1.
\end{equation}
The conditions \eqref{eq:di:per:bnd} reflect the choice \eqref{eq:bck:triplet:p:k}, with the understanding that $\tilde{P}$ is $2L$-periodic. 
Indeed, we resolve function evaluations outside the interval $[0,L]$ by using the identities
\begin{equation}
\label{eq:di:extensions:tild:p}
\tilde{P}
(\tilde{\tau} + 2L) = \tilde{P}(\tilde{\tau}),
\qquad 
\tilde{P}_1(-\tilde{\tau}) = \tilde{P}_1(\tilde{\tau}), 
\qquad
\tilde{P}_2(-\tilde{\tau}) = - \tilde{P}_2(\tilde{\tau}).
\end{equation}
In particular, this means that the triplets \eqref{eq:bck:triplet:p:k} are represented by
\begin{equation}
\label{eq:di:corr:p}
\big(\pbar_1(\xi), \pbar_2(\xi), \omega\big)
\sim \big( \beta_P \tilde{P}_1(L \xi/\pi),
\beta_P\tilde{P}_2(L \xi/\pi),
\pi \omega_P / L \big).
\end{equation}

By picking an appropriate initial condition, we can ensure that $\tilde{P}'_2(0) > 0$ holds whenever $\mu \approx 0$ and $\beta_P$ is sufficiently small. 
Our continuation approach subsequently ensures that we maintain the inequality
\begin{equation}
\nu^{\sigma,\mu}_1 \tilde{P}_1(0)
+\nu^{\sigma,\mu}_2
L\tilde{P}_2'(0)/(\pi \omega_{\sigma,\mu}) > 0
\end{equation}
throughout all our simulations. 
In particular, the sign term in the definition \eqref{eq:per:normalization} for the scaled amplitude $a$ agrees with the sign of $\beta_P$ in view of the correspondence
\eqref{eq:di:corr:p}.
In particular, the numerical equivalent of $a$ is given by the new parameter
\begin{equation}
a \sim \alpha_P := \beta_P
\sqrt{\|\tilde{P}_1\|_\infty^2+ \|\tilde{P}_2\|_\infty^2}.
\end{equation}

\subsection{Diatomic waves}

\begin{figure}[t]
\begin{center}
\includegraphics[width=1\textwidth]{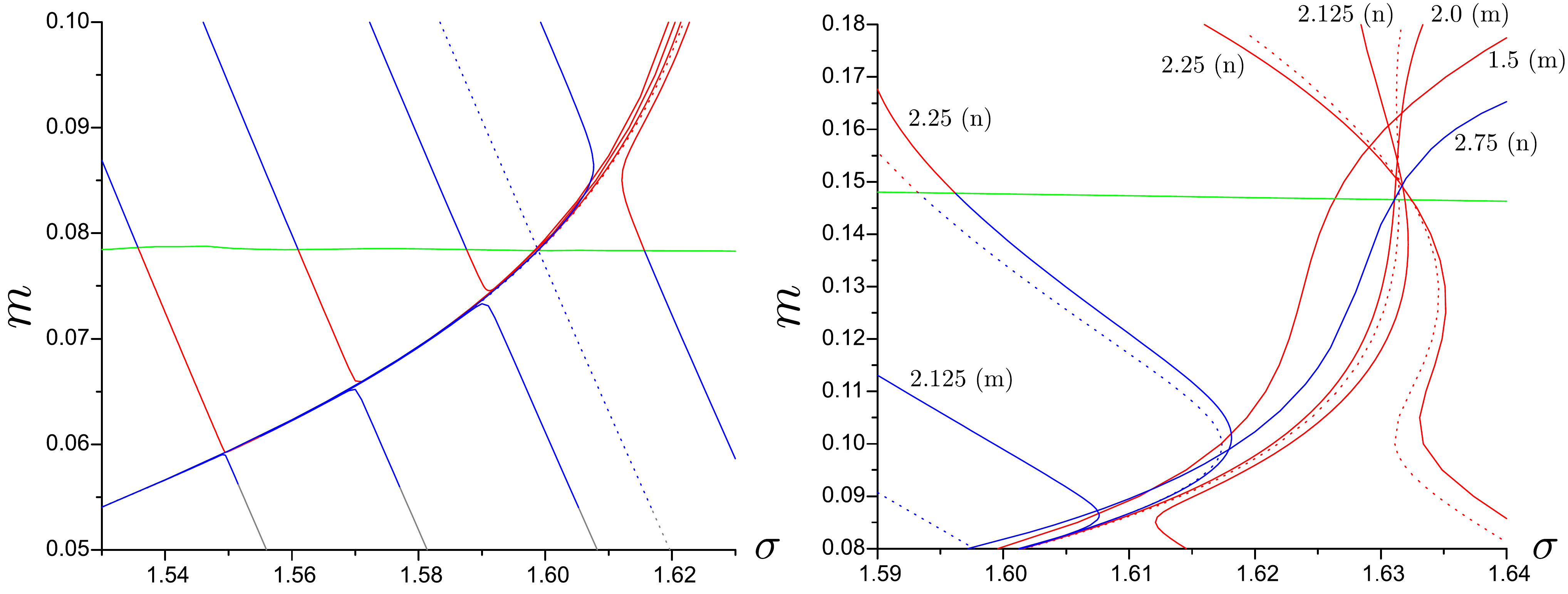}
\caption{Two zoomed views of interesting regions in Fig.\@ \ref{fig:di:overview}. For presentation purposes,
we only plot a subset of the iso-$\kappa$ curves.
The small numbers next to the curves in the right panel correspond to  $\kappa$, while (n) and (m) stand for the nanopteron respectively micropteron branch.
} 
\label{fig:di:overview:zoom}
\end{center}
\end{figure}

Turning to the full ansatz \eqref{eq:bck:beale:ansatz}, we now look for solutions to \eqref{eq:trv;wave:s1:s2} of the form
\begin{equation}
(\sbar_1, \sbar_2)(\xi)
\sim \kappa^2 V(\kappa \xi) + \beta_P P( \kappa \xi).
\end{equation}
Here we have introduced the reparametrization
\begin{equation}
P(\tau) = \tilde{P}(\omega_P \tau / \kappa)
\end{equation}
in order to recast the system in terms of the usual coordinate
\begin{equation}
\tau = \kappa \xi \in [0,L].
\end{equation}
Note that the solitary component $V = (V_1,V_2)$ represents the entire solitary ``core'' of Beale's ansatz \eqref{eq:bck:beale:ansatz}. 
This consists of the ``known'' localized term $(\phi_c,0)$ and the ``unknown error'' term $(w_1,w_2)$, in line with our interpretation of the parameter $\kappa$ in \eqref{kappa approx sol core}. 
In particular, we do not incorporate the monatomic waves from {\S}\ref{sec:mono} because this would only add to the complexity of the numerical procedure.

The full solitary component $V$ should now satisfy the system
\begin{equation}
\begin{array}{lcl}
- \kappa^2 \sigma^2  V'' =
\D_{\kappa}(\mu)
\begin{pmatrix*}
V_1 + \kappa  (V_1^2 + V_2^2)
+ 2 \beta_P( V_1 P_1 + V_2 P_2)
\\[0.2cm]
V_2 + 2 \kappa V_1 V_2
+2 \beta_P (V_1 P_2 + V_2 P_1 )
\end{pmatrix*},
\end{array}
\end{equation}
augmented by the boundary conditions
\begin{equation}
V_1(0) = \frac{1}{8},
\qquad
V_1'(0) = 0,
\qquad
V_2(0) = 0,
\qquad
V_1(L) = 0,
\qquad
V_2(L)  = V_2'(L) = 0.
\end{equation}
The first three of these allow us to extend $V_1$ and $V_2$ as even respectively odd functions, allowing evaluations with $\tau < 0$ to be resolved. 
Evaluations with $\tau > L$ are set to zero, while evaluations of $P$ outside of $[0,L]$ are performed using \eqref{eq:di:extensions:tild:p}.

\subsection{Implementation}
\label{sec:di:impl}
For each individual run of the collocation solver, we fix the parameter $\kappa > 0$ together with one of the variables from the set $\{\sigma , \mu, \beta_P\}$.
The remaining two variables then need to be computed, along with $\omega_P$ and the functions $(V_1, V_2, P_1, P_2)$.
Since the latter all satisfy second-order MFDEs, our solver requires $4\times2 + 3 = 11$ boundary conditions, which corresponds with the eleven boundary conditions formulated above.
This freedom to choose the second fixed parameter is essential for our continuation approach, since it allows us to move past fold points by switching our choice.

In Fig.\@ \ref{fig:di:overview} we consider several values of $\kappa \in \mathbb{N}/8$ and trace out the corresponding curve(s) in the $(\sigma,m)$ landscape where we were able to find diatomic solitary waves for the chosen value of $\kappa$. 
In addition, we use broken lines to plot these curves for the special values 
\begin{equation}
\label{eq:di:spec:kappa}
\kappa \in \{ 2.0515,2.237567  \}
\end{equation}
where the behavior of these curves experiences a structural change.
These values were determined by an unsophisticated bisection approach.
In Fig.\@ \ref{fig:di:overview:zoom} we provide a zoomed-in view of two interesting areas near the folds, which we discuss in further detail below.

\begin{figure}[t]
\begin{center}
\includegraphics[width=1\textwidth]{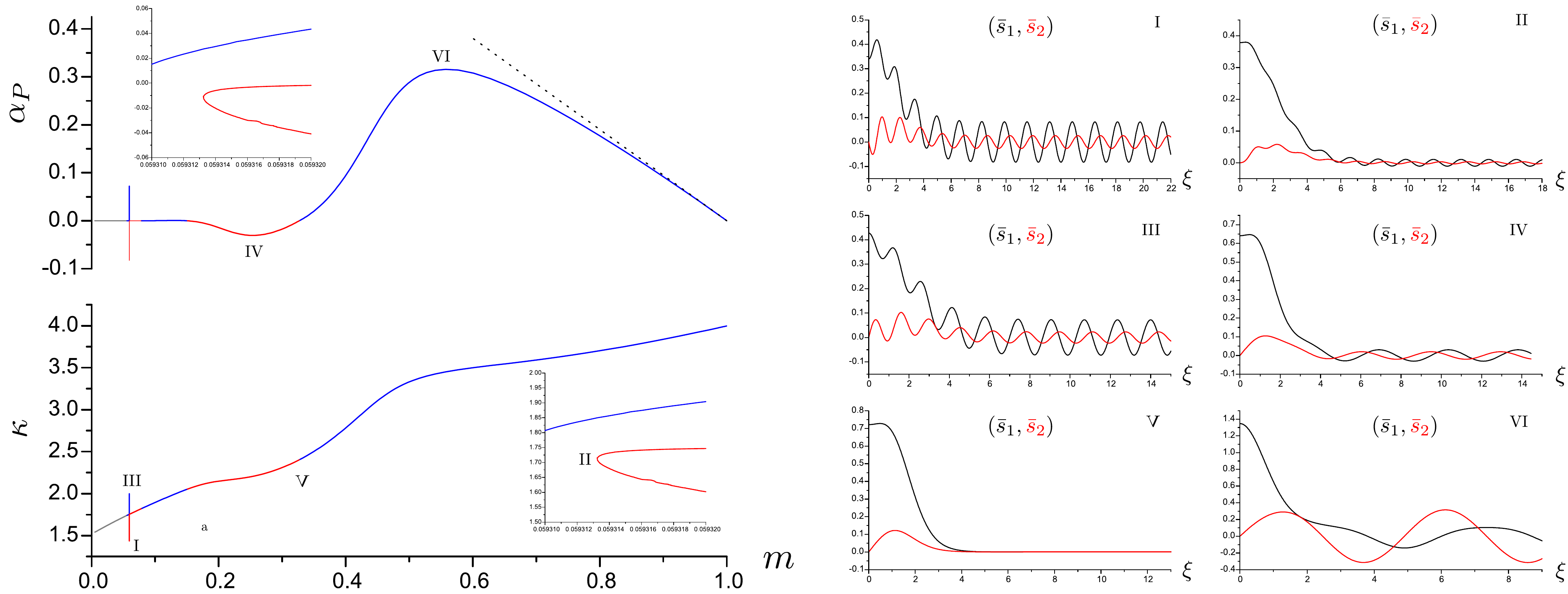}
\caption{
The left graphs contain vertical
cross-sections of Fig.\@ \ref{fig:di:overview} 
with $\sigma = 1.55$. The dashed line represents
the leading-order ripple-amplitude prediction
\eqref{eq:bck:expr:coeff:ampl} computed with
\eqref{eq:mono:ampl:coeff}, while the insets
zoom in on the region where the micropteron
and nanopteron branches approach each other.
The six graphs on the right contain the
solution profiles associated to the special
points I through VI marked on the left graphs.
} 
\label{fig:cross:c:1p55}
\end{center}
\end{figure}

\begin{figure}[t]
\begin{center}
\includegraphics[width=1\textwidth]{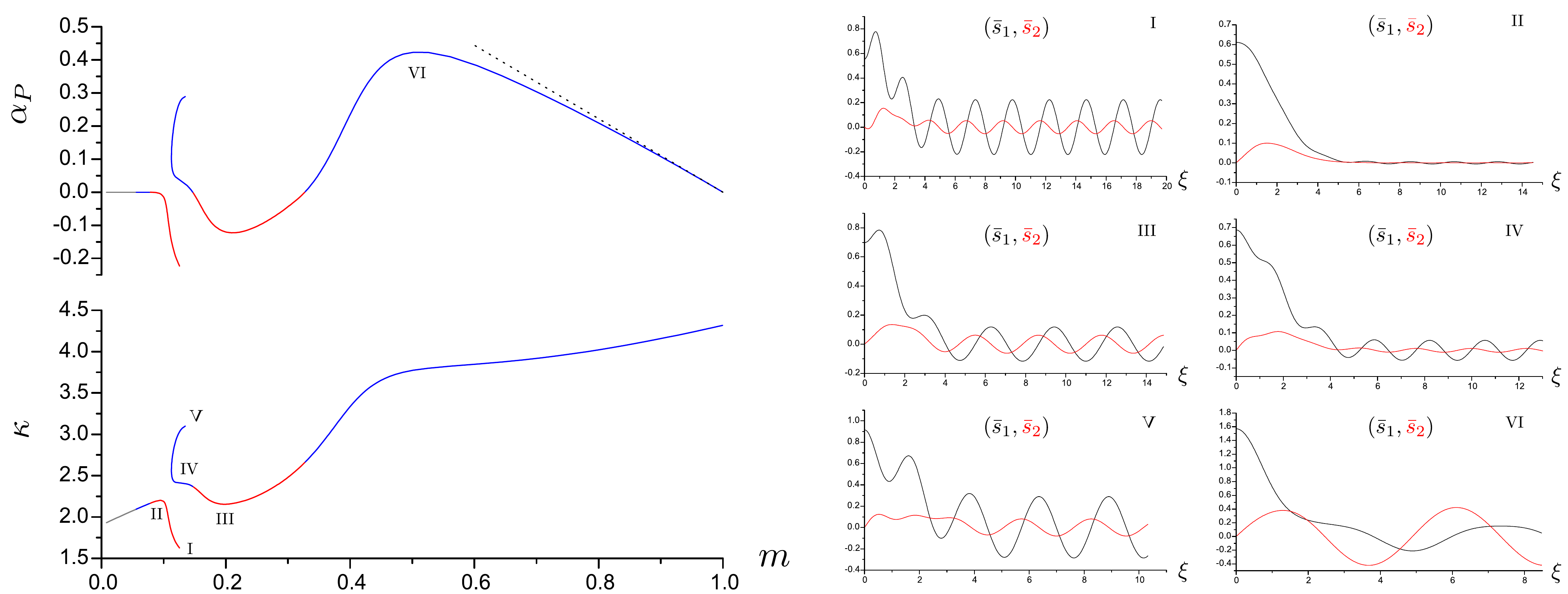}
\caption{Cross-sections and solution profiles at $\sigma = 1.625$; see Fig.\@ \ref{fig:cross:c:1p55} for more information.
} 
\label{fig:cross:c:1p625}
\end{center}
\end{figure}

\begin{figure}[t]
\begin{center}
\includegraphics[width=1\textwidth]{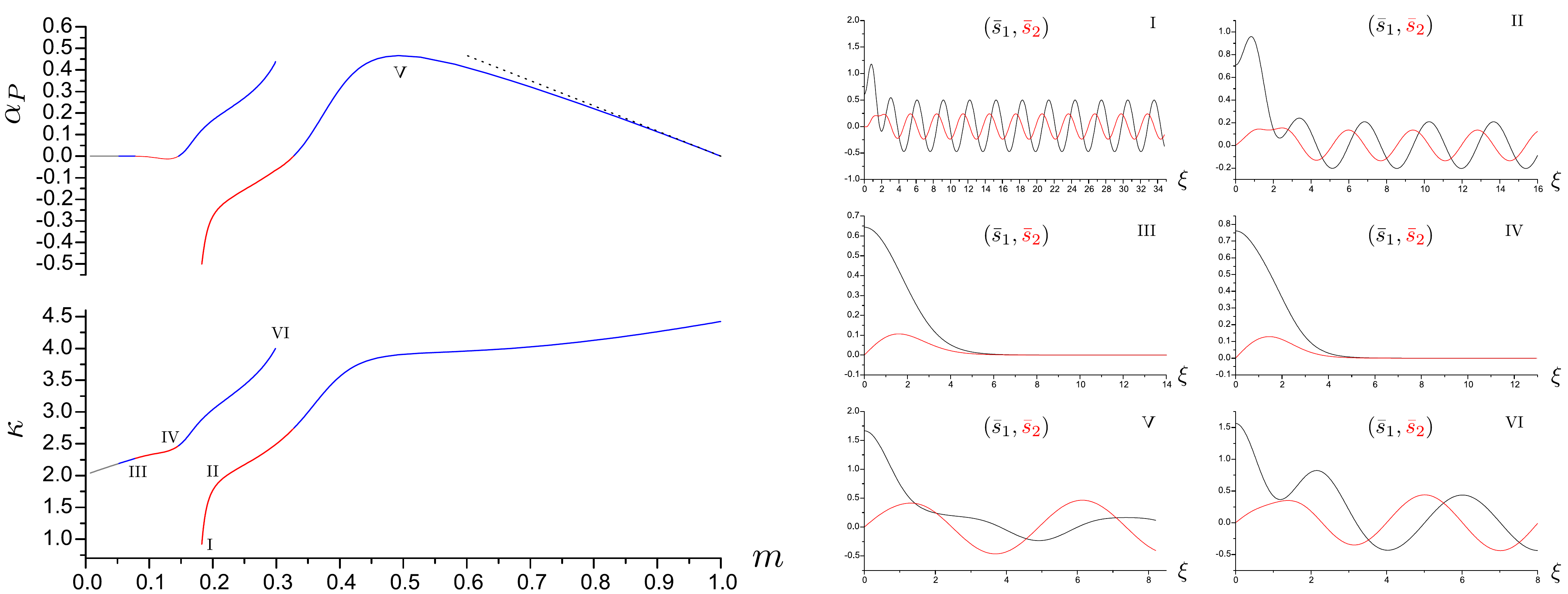}
\caption{Cross-sections and solution profiles at $\sigma = 1.65$; see Fig.\@ \ref{fig:cross:c:1p55} for more information.
} 
\label{fig:cross:c:1p65}
\end{center}
\end{figure}

We use the ripple-amplitude $\alpha_P$ to color the curves. 
In principle, we use blue and red for positive respectively negatives values. 
However, we pay special attention to the regime where the ripple amplitude is at least a factor of $10^{-5}$ smaller than the center amplitude of the solitary component, i.e., where
\begin{equation}
\label{eq:di:small:alpha:p}
|\alpha_P| < 10^{-5} \kappa^2 V_1(0)  =  10^{-5} \frac{\kappa^2}{8}.
\end{equation}
In particular, whenever \eqref{eq:di:small:alpha:p} holds for an interval of $\mu$ of length at least $0.01$, we color the entire segment of the curve where it holds gray. 
In this so-called ``small-ripple'' regime it is hard to distinguish numerically between positive and negative values of $\alpha_P$, as can be seen from the top-left graphs in Figs.\@ \ref{fig:cross:c:1p55}--\ref{fig:cross:c:1p65}.

This problem is further illustrated by the green curves in Figs.\@ \ref{fig:di:overview} and \ref{fig:di:overview:zoom}, which we computed by fixing $\alpha_P = \beta_P = 0$ and performing a $\kappa$ scan. 
In particular, the diatomic waves along these curves are in fact solitary waves.
We were only able to continue these branches slightly into the ``small-ripple'' regime, after which the scheme failed to converge.
We emphasize that this does not necessarily mean that these branches terminate.

The bottom-left graphs in Figs.\@ \ref{fig:cross:c:1p55}--\ref{fig:cross:c:1p65} can be seen as vertical cross-sections of Fig.\@ \ref{fig:di:overview}. 
Indeed, they were obtained by fixing the wave speed $\sigma$ and performing (several) $\kappa$-scans, using continuation to find appropriate initial solutions.
Besides $\kappa$, these figures also visualize $\alpha_P$ as a function of $m$. 

As can be seen, in these cross-sections we were unable to access regions with $\kappa \le 1$. 
We suspect that this is a consequence of the increasing number of mesh intervals that are required to resolve the oscillations in the solitary component $V$, which leads to memory issues. 
Indeed, at present our software uses legacy 32-bit code that limits the amount of accessible memory to roughly 4Gb. 
Due to the (relatively) large number of components in the system and the non-standard structure of the underlying matrices caused by the shifted arguments in the MFDE \eqref{tw eqns}, this limit is reached much sooner than one would encounter when solving ODEs.

For similar reasons, we have not been able to extend the (blue) fold curve that presently terminates at $(\sigma,m) = (1.46518,0.038711)$ in Fig.\@ \ref{fig:di:overview} to smaller values of $\sigma$. 
More precisely, we could not continue the $\kappa = 1.5$ ``nanopteron'' branch (emanating from $m=0$) past this point. 
On the other hand, for smaller values of $\kappa$ we were unable to resolve the ``turn'' where the ``micropteron'' (emanating from $m=1$) and ``nanopteron'' branches move away from each other.
Indeed, Fig.\@ \ref{fig:di:overview:zoom}a clearly shows that these ``turn regions'' become increasingly thin as $\kappa$ decreases. Increasingly delicate techniques are therefore required to prevent the
continuation procedure from simply jumping between the two branches. 

\subsection{Discussion}
The overview in Fig.\@ \ref{fig:di:overview} and the iso-$\sigma$ cross-sections in Figs.\@ \ref{fig:cross:c:1p55}--\ref{fig:cross:c:1p65} clearly show that one cannot simply speak about
separate  ``micropteron'' and ``nanopteron'' surfaces.
For example, the iso-$\kappa$ curves converge precisely to the critical speed-of-sound $c^{(s)}_m$ defined in \eqref{lw speed of sound} for the long-wave nanopterons as $\kappa \downarrow 0$.
In addition, the curves emanating from $m=1$ and $m=0$ are connected together in a complicated fashion. 

However, if one fixes either of the parameters $\kappa$ or $\sigma$, the resulting one-dimensional cross-sections do generically consist of two separate curves. 
We have identified several exceptions in Figs.\@ \ref{fig:di:overview} and \ref{fig:di:overview:zoom} and conjecture that there is a countable set of such special values that accumulates at  $\kappa = 0$ respectively $\sigma = \sqrt{2}$.

In order to support this conjecture, let us consider the region where $\sigma \le 1.61$ and $m \le 0.08$, where a large set of iso-$\kappa$ curves converge together to form a narrow ``fold-region'' that becomes increasingly thin as $\sigma$ decreases.
Depending on the sign of $\alpha_P$, the ``micropteron'' and ``nanopteron'' curves seem to bend sharply to the left or right as this fold is approached. 
Our observations for $\kappa = 2.0515$ suggest that a switch in direction occurs precisely when the ``micropteron'' and ``nanopteron'' curves connect to each other and the associated $\alpha_P$ function admits a double root when passing through the ``connection'' point at $m \sim 0.0784$. We found faint numerical hints of a second such direction switch when the fold crosses through  $m \sim 0.05$, but did not have the numerical resolution to fully resolve this bifurcation.

In particular, we suspect that the (green) branches of solitary waves play an important role as they contain the ``crossing-points'' between the nanopteron and micropteron subsurfaces. 
Although we had trouble tracking them deep into the ``small-ripple'' regime, we did find several of these branches. 
They indeed appear to accumulate in the small-mass regime, lending credence to the conjectures discussed above.
It interesting to note that these branches are not horizontal, i.e., the value of $m$ changes gradually as $\sigma$ is varied. 
This is especially true for the $m \sim 0.15$ branch, which has a relatively substantial slope.

We believe that the fold-region discussed above continues into the corner $(\sigma, m) = (\sqrt{2}, 0)$, admitting a countable number of connection switches as it crosses through the branches of solitary waves.
In our opinion it would be extremely interesting to perform a theoretical analysis near this corner.
One could proceed by examining how the long-wave techniques developed in \cite{faver-wright} and the small-mass argument discussed in \cite{hoffman-wright} break down as their auxiliary parameters approach the limits $1/m \to \infty$ respectively $\sigma \to \sqrt{2}$. 
Combining these approaches with the insights developed by Lombardi \cite{lombardi} to uncover exponentially small phenomena could hopefully lead to some useful insights here.

A second (related) question  that deserves further attention concerns the behavior of the ``micropteron'' and ``nanopteron'' branches after their ``near-collision'' events; see the left graphs in Figs.\@ \ref{fig:cross:c:1p55}--\ref{fig:cross:c:1p65}. It appears that the red branches where $\alpha_P < 0$ suffer a collapse in the amplitude of the solitary component of the wave, i.e., $\kappa \downarrow 0$, potentially converging to a branch of purely periodic ripples. 
It is unclear what happens to the blue branches, where both the ripple amplitude and the core amplitude experience significant growth.

\section{Stability of Solitary Waves}
\label{sec:full:fpu}

\begin{figure}[t]
\begin{center}
\includegraphics[width=1\textwidth]{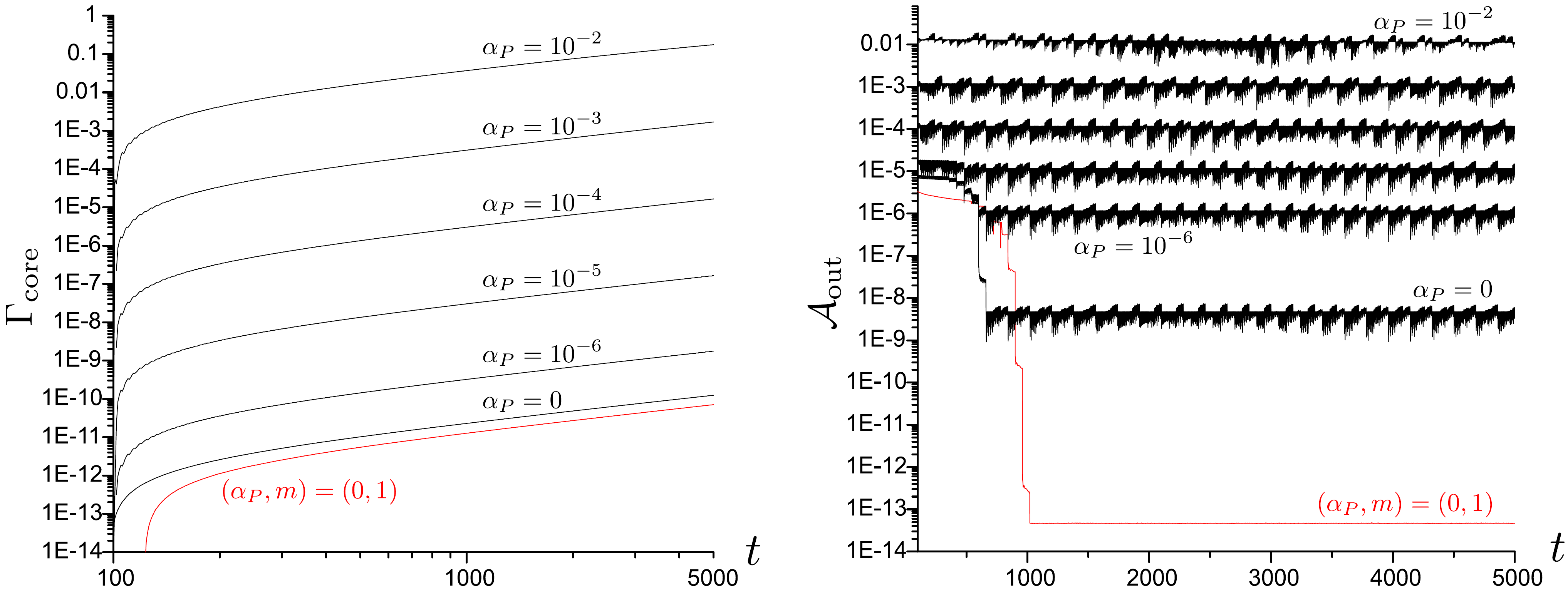}
\caption{behavior of the core energy loss 
\eqref{eq:fput:core:energy}
and the outer amplitude \eqref{eq:fput:outer:ampl} for several simulations of the full FPUT system \eqref{eq:fpu:r:p:coords}. The values of $\alpha_P$ correspond with the initial conditions
\eqref{eq:fput:init:conds}. The red curves represent the monatomic equal-mass solitary wave at $\kappa = 5/2$.
} 
\label{fig:full:fput}
\end{center}
\end{figure}

In order to gain some insight into the stability of the solitary waves that we found in {\S}\ref{sec:di}, we performed a series of direct time-integrations of the full FPUT system. 
For our purposes here,it is advantageous to introduce a new momentum variable $p_j = \dot{x}_j$ and reformulate \eqref{newton} as 
\begin{equation}
\label{eq:fpu:r:p:coords}
\dot{r_j} = p_{j+1} - p_j,
\qquad \qquad m_j \dot{p}_j = F( r_j) -
F(r_{j-1}),
\end{equation}
recalling that the masses are given by \eqref{masses}.
The advantage of these coordinates is that the energy function
\begin{equation}
\mathcal{E}_{\Lambda}
= \sum_{j \in \Lambda} 
\Big[
\frac{1}{2} m_j p_j^2
+\frac{1}{2} r_j^2 + \frac{1}{3} r_j^3
\Big]
\end{equation}
is conserved in time upon taking
$\Lambda = \mathbb{Z}$,
but also easy to compute for any subset $\Lambda \subset \mathbb{Z}$.
In general, diatomic waves  of the form \eqref{eq:bck:beale:ansatz} will have infinite energy on account of the periodic background state. 
However, for solitary waves (where $a =0$) the energy is finite, since the pair $(w_1,w_2)$ is exponentially localized.

In the monatomic case $m=1$, Friesecke and Pego \cite{friesecke-pego2,friesecke-pego3,friesecke-pego4} established that the traveling wave solution $r_j(t) = \phi_{\cep}( j - \cep t)$ is stable under the dynamics of \eqref{eq:fpu:r:p:coords} whenever $\ep > 0$ is sufficiently small. In particular, any ``sufficiently small'' initial perturbation to the wave will die out over time, although the speed and phase of the wave could be slightly changed.
In practice, one sees that perturbations produce a small wrinkle that separates itself from the core of the wave and travels at a slower speed.

By constrast, the stability properties of diatomic traveling waves are poorly understood at present. 
The main obstruction is that one has to control the fluctuations caused by the ripples; see, e.g., \cite{johnson-wright} for a discussion of the complications that arise in the stability analysis of water wave nanopterons. 
Indeed, the numerical results from \cite{co-ops} --- which use the suitably scaled solitary KdV profiles \eqref{eq:int:def:Phi:m} as initial conditions for \eqref{eq:fpu:r:p:coords} ---  clearly indicate that ripples will appear that \textit{do not} detach from the core of the wave.
Instead, they slowly drain energy from this core, causing the amplitude to decay over time. 
Nevertheless, these structures persist over much longer timescales than those that can be extracted by using amplitude  equations. 
For example, Gaison, Moskow, Wright and Zhang \cite{gmwz} consider a general class of polyatomic lattices that includes our case here and show that the relevant KdV reductions that govern the long-wave limit remain valid over algebraically long timescales.

In order to examine these issues, we introduce the time-dependent set of gridpoints
\begin{equation}
\Lambda_{\mathrm{core}}(t) = \set{j \in \Z}{ \big|j - \mathrm{argmax} |r_j(t)| \big| \le 20 } 
\end{equation}
that is centered around the peak of a solution to \eqref{eq:fpu:r:p:coords}. 
Our goal is to monitor the behavior of the loss function
\begin{equation}
\label{eq:fput:core:energy}
\Gamma_{\mathrm{core}}(t) = \frac{
\mathcal{E}_{\Lambda_{\mathrm{core}}(100)}
-
\mathcal{E}_{\Lambda_{\mathrm{core}}(t)}
}{
\mathcal{E}_{\Lambda_{\mathrm{core}}(100)}
}
\end{equation}
for various solutions resembling diatomic and monatomic waves.
This fraction measures the relative energy loss in the core of the wave compared to the situation at $t = 100$, allowing sufficient time for initial transients to decay.
We deliberately do not use the peak amplitude of the wave here, since at each point in time only a discrete subset of the underlying smooth waveprofile is ``sampled'' on the lattice.
We do however keep track of the outer amplitude
\begin{equation}
\label{eq:fput:outer:ampl}
\mathcal{A}_{\mathrm{out}}(t) = \max_{j \notin \Lambda_{\mathrm{core}}(t)} |r_j(t)|
\end{equation}
as a secondary measure.

\subsection{Implementation}

We use the \texttt{solve{\_}ivp} method from the SciPy package for Python to integrate the problem \eqref{eq:fpu:r:p:coords} on the grid $\Lambda_{\mathrm{grid}} = \{1, 2, \ldots, 400 \}$.
In particular, we set both $r$ and $p$ to zero outside this grid. 
We consider six separate initial conditions, which each consist of the solitary part (i.e., the pair $(V_1, V_2)$) of one of the diatomic waves computed in {\S}\ref{sec:di}. 
More precisely, we pick $\kappa = 5/2$ and consider the waves associated to the pairs
\begin{equation}
\label{eq:fput:init:conds}
\begin{array}{lcl}
(\alpha_P,m)
&\in & \{ (10^{-2},0.33797458), (10^{-3},0.32800968),
(10^{-4},0.32711659),
\\[0.2cm]
& & \qquad
(10^{-5},0.32702829),
(10^{-6},0.32701947),
(0,0.32701849) \},
\end{array}
\end{equation}
which have wave speeds in the range $\sigma \in [1.563,1.578]$.
For comparison purposes, we also considered the monatomic wave $(\alpha_P,m) = (0, 1)$, again with $\kappa = 5/2$.
The results of these simulations can be found in Fig.\@ \ref{fig:full:fput}, where we plot the evolution of the core energy loss $\Gamma_{\mathrm{core}}(t)$ and the outer amplitude $\mathcal{A}_{\mathrm{out}}(t)$.

We note that the \texttt{solve{\_}ivp} routine utilizes the RK4 scheme with an adaptive step-size $\Delta t$.
We monitored the energy $\mathcal{E}_{\Lambda_{\mathrm{grid}}}$ over the full grid in order to keep track of potential discretization errors introduced by the scheme. 
This turned out to be a crucial precaution, because the step-size automatically chosen  by \texttt{solve{\_}ivp} led to unacceptable fluctuations in this energy. 
To prevent this, we manually enforced the step-size restriction $\Delta t \le 10^{-3}$, which is considerably lower than in \cite{co-ops}. 
We note that an alternative approach could be to use symplectic energy-preserving schemes as described in \cite{hairer2006geometric}, but we believe that this is too cumbersome for our illustrative purposes here.

We run our simulations for $t \in [0, 5000]$, which means that the waves will have shifted roughly 7850 lattice points to the left. 
Since this considerably exceeds the size of our grid, we need to account for this movement in a special fashion.
Our choice here is to use a ``windowing'' procedure, where we shift the solution rightwards to recenter the peak at the center of the computation grid, filling $r$ and $p$ with zeros at the empty positions on the left side of the grid. 
In order to prevent large discontinuities arising from the sudden cutoff at the right end of the grid, we apply the pointwise multiplication
\begin{equation}
(r_i,p_i) \mapsto   e^{-y_i^2/(1 - y_i^2)} (r_i, p_i), \qquad \qquad
y_i = \max\{ (i - 300)/100 , 0\}.
\end{equation}
In particular, we use a smooth  cut-off function to gradually scale the solution on the right $25\%$ of the lattice sites down to zero.
In contrast to the approach in \cite{co-ops}, we only perform this recentering and windowing procedure once per sixty units of time. 
We emphasize that other methods are available to deal with this problem, such as the freezing technique developed by Beyn and coworkers \cite{beyn2004freezing}.

\subsection{Discussion}
Naturally, due to discretization effects and rounding issues there will always be some energy leakage as a numerical wave moves through the lattice. 
In addition, sampling and interpolation errors occur when passing initial conditions from the boundary value problem solver used in {\S}\ref{sec:di} to the FPUT simulator discussed here. 
This causes the initial transient behavior and subsequent slow energy leakage that is displayed in Fig.\@ \ref{fig:full:fput} for the monatomic wave.
We use this as a baseline to interpret the behavior of the diatomic waves.

The results in Fig.\@ \ref{fig:full:fput}a show that our numerical diatomic solitary wave loses energy at an extremely slow rate that is comparable to its monatomic counterpart. 
In addition, the initial transient clearly indicates that our diatomic solitary wave is stable in a certain sense.  
Notice furthermore that this coherence disappears rapidly if the mass $m$ is disturbed. 
Indeed, the size of the ripple-amplitude $\alpha_P$ is clearly correlated with the speed at which the core of the wave loses energy.

We reiterate that the results in \cite{co-ops} already suggest that --- in general --- diatomic waves decay at rates that are much slower that those suggested by their formal KdV approximations. 
This effect is amplified in the small-amplitude regime, where \cite{co-ops} contains examples that display practically no decay on very long time-scales. 
For this reason, we consider the relatively large value $\kappa = 5/2$, which through the parameter translation 
\begin{equation}
(\kappa,m) \sim (\sqrt{24} \tilde{\ep}, 1/m_2)
\end{equation}
means that our results should be compared to the results in \cite{co-ops} with $\tilde{\ep} =1/2$ and $m_2 = \pi$, which were the largest values that the authors consider. 
In particular, our results should be compared with the top-left plot in \cite[Fig.\@ 8]{co-ops}, noting that our time-interval corresponds roughly with $[0, 10^4]$ in that figure. 

The conclusion from this comparison is that the  observed ``outer amplitude'' $\mathcal{A}_{\mathrm{out}}$ for our diatomic solitary waves is orders of magnitude smaller than the wakes observed in \cite{co-ops}. 
The loss function $\Gamma_{\mathrm{core}}$ exhibits a similar scale reduction, although this is harder to read-off from the figure due to the smaller numbers. 
Together, we feel that these observations strongly suggest that the $m \sim 0.33$ branch of solitary waves found in {\S}\ref{sec:di} are indeed solitary and stable, providing a robust transport mechanism in the diatomic setting.

\section{Future Directions}
\label{sec:disc}

In this paper we considered the numerical behavior of diatomic FPUT lattices with the quadratic spring force $F(r) = r+r^2$.
The MiM lattice is a natural candidate for future investigations of this sort, since its traveling wave equations are MFDEs similar to \eqref{tw eqns}, and since it possesses a sturdy theory of solitary waves \cite{ksx, faver-goodman-wright} and nanopterons \cite{faver-mim-nanopteron} in the small bead-resonator mass limit.
In particular, from \cite{ksx} and \cite{faver-goodman-wright} there is a simple, explicit formula for those mass ratios accumulating at 0 at which the MiM lattice has solitary waves.
It is also possible to pose a ``stiff internal spring'' limit for the MiM lattice, in which the bead-resonator spring becomes arbitrarily stiff, and the MiM lattice again reduces to a monatomic FPUT lattice.
Solitary waves are known to exist in this limit, too \cite{faver-goodman-wright}, and one expects nanopterons there was well.
There could be an interesting parameter overlap, similar to the regimes in Fig.\@ \ref{fig-all_nonlocals}, in which the bead-resonator spring force is extremely stiff and the bead-resonator mass ratio is extremely small.

While we numerically simulated periodic solutions to FPUT lattices, we did not consider how different families of periodics relate to each other, in the spirit of the full nanopterons and micropterons in Fig.\@ \ref{fig-all_nonlocals}.
Friesecke and Mikikits-Leitner \cite{fml} construct periodic traveling waves in the long wave limit for monatomic lattices.  
Do these have any connection to the periodic waves in diatomic lattices from \cite{faver-wright} when the mass ratio is close to 0 or 1?
We also mention that Betti and Pelinovsky study periodic waves in a diatomic lattice with Hertzian spring forces \cite{betti-pelinovsky}.
They begin with periodic solutions in the small mass limit but manage to extend them numerically to the equal mass limit.
For a given wave speed $c$, is it possible to extend the FPUT periodics in the same way, from $m \approx 0$ to $m \approx 1$?

Last, in more general polyatomic FPUT lattices, in which both the masses and the spring forces repeat with some finite periodicity, it is known that solutions to the equations of motion with suitably scaled initial data look like KdV $\sech^2$-type solitary waves over long times \cite{gmwz}.
However, it is not yet known whether these solitary wave approximations persist for all time, as in the the monatomic lattice, evolve into nanopterons, as in the diatomic lattice, or become something else entirely.
In a more complicated polyatomic lattice, there are fewer opportunities for a natural ``material'' limit, like the small or equal mass regimes, to reduce the polyatomic lattice to monatomic.
Nonetheless, our numerical methods could give insight into the formation of, at least, solitary waves in the long wave polyatomic limit.

\section*{Acknowledgments}

Funding: Both authors acknowledge support from the Netherlands Organization for Scientific Research (NWO) (grant 639.032.612).

\bibliographystyle{siam}
\bibliography{numerics_bib}

\end{document}